\documentclass[final]{siamltex}

\usepackage[normalem]{ulem}
\usepackage{amsmath}
\usepackage{amsfonts}
\usepackage{amssymb}
\usepackage{graphicx}
\usepackage{epstopdf,euscript}
\usepackage{url}
\usepackage{subfigure}
\usepackage{multirow}
\usepackage{array}
\usepackage{color}

\def\bq{\begin{quotation}}
\def\eq{\end{quotation}}

\def\D{\Delta}

\def\h{\eta}

\def\V{\Theta}

\def\m{\mu}
\def\n{\nu}

\def\o{\omega}

\newcommand{\bfDe}{{\boldmath \D}}
\newcommand{\bbM}{{ \mathbb{M} }}

\def\2nm#1{\|#1\|_2}

\def\Ra#1{\mathrm{Range}(#1)}

\newcommand{\ars}[1]{\left[ \begin{array}{#1}}
\newcommand{\are}{\end{array} \right] }
\newcommand{\oars}[1]{\begin{array}{#1}}
\newcommand{\oare}{\end{array}}
\newcommand{\rars}[1]{\left( \begin{array}{#1}}
\newcommand{\rare}{\end{array} \right) }

\newcommand{\eqs}{\begin{eqnarray}}
\newcommand{\eqe}{\end{eqnarray}}
\newcommand{\eqsn}{\begin{eqnarray*}}
\newcommand{\eqen}{\end{eqnarray*}}

\newcommand{\bmp}[2]{\begin{minipage}#1{#2}}
\newcommand{\emp}{\end{minipage}}

\newcommand{\ens}{\begin{enumerate}}
\newcommand{\ene}{\end{enumerate}}

\newcommand{\its}{\begin{itemize}}
\newcommand{\ite}{\end{itemize}}

\newcommand{\des}{\begin{description}}
\newcommand{\dee}{\end{description}}

\def\defs{\begin{definition}}
\def\defe{\end{definition}}
\def\teos{\begin{theorem}}
\def\teoe{\end{theorem}}
\def\prfs{\begin{proof}}
\def\prfe{\end{proof}}
\def\exas{\begin{exampl}}
\def\exae{\end{exampl}}
\def\excs{\begin{exercise}}
\def\exce{\end{exercise}}
\def\cors{\begin{corollary}}
\def\core{\end{corollary}}

\def\wh{\widehat}

\newcommand{\Hinf}{{\mathcal{H}_{\infty}}}

\newcommand{\cbfm}{\mbox{\boldmath${\mathit{m}}$} }

\newcommand{\bfsfp}{\mbox{\boldmath$\mathsf{p}$} }

\newcommand{\comment}[1]{} 

\newcommand{\bfPsi}{\boldsymbol{\Psi}}

\newcommand{\bfpi}{\boldsymbol{\pi}}

\newcommand{\bfchi}{\mbox{\boldmath$\chi$}}

\newcommand{\bfp}{{\bf p}}
\newcommand{\bfx}{{\bf x}}

\newcommand{\bfy}{{\bf y}}

\newcommand{\bfu}{{\bf u}}

\newcommand{\bfv}{{\bf v}}

\newcommand{\bfA}{{\bf A}}
\newcommand{\bfB}{{\bf B}}
\newcommand{\bfC}{{\bf C}}
\newcommand{\bfD}{{\bf D}}
\newcommand{\bfE}{{\bf E}}

\newcommand{\bfI}{{\bf I}}

\newcommand{\bfY}{{\bf Y}}
\newcommand{\bfU}{{\bf U}}
\newcommand{\bfM}{{\bf M}}
\newcommand{\bfZ}{{\bf Z}}
\newcommand{\bfV}{{\bf V}}
\newcommand{\bfW}{{\bf W}}

\newcommand{\IR}{{\mathbb{R}}}
\newcommand{\IC}{{\mathbb{C}}}

\newcommand{\bea}{\left[ \begin{array} }
\newcommand{\eea}{ \end{array} \right] }

\newcommand{\data}{ {\mathbb{D}} }

\title{Nonlinear Parametric Inversion using Interpolatory Model Reduction\footnotemark[1]}

\author{Eric de Sturler\footnotemark[2]\ ,
Serkan Gugercin\footnotemark[2]\ ,
Misha E. Kilmer\footnotemark[3]\ ,
Saifon Chaturantabut\footnotemark[2]\ ,
Christopher Beattie\footnotemark[2]\ ,
and Meghan O'Connell\footnotemark[3]}

\begin{document}
\maketitle
\renewcommand{\thefootnote}{\fnsymbol{footnote}}
\footnotetext[1]{This material is based upon work supported by the National Science
Foundation under Grants No. {NSF-DMS} 1025327, {NSF DMS} 1217156 and 1217161,
{NSF-DMS} 0645347, and {NIH} R01-CA154774.}
\footnotetext[2]{Department of Mathematics, Virginia Tech, Blacksburg, VA 24061.}
\footnotetext[3]{Department of Mathematics, Tufts University, Medford, MA 02115.}
\renewcommand{\thefootnote}{\arabic{footnote}}

\begin{abstract}
Nonlinear parametric inverse problems appear in several prominent
applications; one such application
is Diffuse Optical Tomography (DOT) in medical image reconstruction.
Such inverse problems present huge computational challenges,
mostly due to the need for solving a sequence of large-scale discretized,  
parametrized, partial differential equations (PDEs) in the forward model.
In this paper, we show how interpolatory parametric model reduction  can significantly
reduce the cost of the inversion process in DOT by drastically reducing the
computational cost of solving the forward problems.
The key observation is that
function evaluations for the underlying optimization problem
may be viewed as transfer function evaluations along the imaginary axis; a similar observation
holds for  Jacobian evaluations as well. This motivates the use of
system-theoretic model order reduction methods. We discuss the construction
and use of interpolatory parametric reduced models as surrogates for the full
forward model. Within the DOT setting, these surrogate models can approximate
both the cost functional and the associated Jacobian with very little loss of
accuracy while significantly reducing the cost of the overall inversion process.
Four numerical examples
illustrate the efficiency of the proposed approach.
Although we focus on DOT in this paper, we believe that our approach
is applicable much more generally.

\end{abstract}



\pagestyle{myheadings}
\thispagestyle{plain}
\markboth{DE STURLER, GUGERCIN, KILMER, CHATURANTABUT, BEATTIE, O'CONNELL}
{PARAMETRIC INVERSION USING MODEL REDUCTION}


\section{Introduction}\label{sec:intro}


Nonlinear inverse problems, as exemplified by medical image reconstruction 
or identification and localization 
of anomalous regions (e.g., tumors in the body \cite{bushberg2003essential},
contaminant pools in the earth \cite{james33optimal},
or cracks in a material sample \cite{stavroulakis2001inverse}), are 
commonly encountered yet
remain very expensive to solve.
In such inverse problems, 
one wishes to recover information identifying an unknown spatial distribution 
(the \emph{image}) of some quantity of interest 
within a given medium that is not directly observable. 
For example, identifying anomalous regions of electrical conductivity 
in a sample of muscle tissue aids in identification and localization of potential tumor sites. 
  
The principal tool linking the images that are sought to correlated data that may be observed and measured 
 is a mathematical model, the {\it forward model}. 
Within the context considered here,  
forward models are large-scale, discretized, 2D or 3D, partial differential equations.
These forward models constitute the functions to be evaluated 
for the underlying optimization problem that fits images of interest to observed data, 
so it is necessary to resolve 
these large-scale forward problems many times in order to 
 to recover and reconstruct an image to some desired resolution. 
This constitutes the largest single
computational impediment to effective, practical use of some imaging
modalities and low quality image resolution is an all too common and regrettable outcome. 
Rapid advances in technology make it possible to take
many more measurements, which allows for higher resolution reconstructions in principle.
Yet, the advantage of these additional measurements may not be realized in practice, since solving the forward
problem at compatible resolutions may remain formidably expensive.
The features outlined above apply to many
inverse problems found in fields such as geophysics
\cite{zhdanov2002geophysical,snieder1999inverse}, medical imaging
\cite{louis1992medical,webb2003introduction,arridge}, hydrology
\cite{yeh22review,carrera2005inverse,sun1994inverse}, and nondestructive
evaluation \cite{marklein2002linear,liu2003computational}. 
The need for
innovative, efficient, and accurate algorithms for such problems is as great as
ever. 

In this paper, we propose to combine effective low-order parametric image representations
with techniques developed for
interpolatory model reduction  in order to reduce drastically the
computational cost of solving the nonlinear inverse problem.
We assume that the fundamental objective is the
identification and characterization of anomalous regions in an otherwise
nearly homogeneous medium and that observations may be made  
only at the boundary of the medium.
For such inverse problems,
low-order parametric
image models are able to capture both the edge geometry of the 
anomalies as well as the intensities of quantities of interest
within the anomalies;
 one need only 
recover a relatively small number (compared to the number of grid points)
of parameter values defining these edges
and intensities.
We adopt the Parametric Level Set (PaLS) approach
described in \cite{Aghasi_etal11}.
In this approach, the inverse problem becomes one of finding a
parameter vector $\bfsfp$ that satisfies
\eqs \label{eq:GenInvProbl}
\bfsfp := \arg \min_{\bfsfp \in \mathbb{R}^\ell} \| \bbM(\bfsfp) - \data \|_2 ,
\eqe
where $\bbM(\bfsfp)$ denotes the synthetically generated data based on
a (regularized) forward model for a given input parameter vector
$\bfsfp$, and $\data$ is the data vector comprised of
the measurements at the detectors. In our application, the
measurements are taken in the frequency domain.
The optimization problem is typically
solved by some nonlinear least squares method, in particular,
we will use the TREGS algorithm \cite{StuKil11c}.
Although the inverse problem in this parametrized framework is
considerably easier to solve than solving for the values of the unknown function at
every grid point in a 2D or 3D grid, the forward solves required for function
and Jacobian evaluations still lead to significant computational effort.
In this paper, 
model reduction is employed to bring down the cost of forward
solves, making the overall inversion process considerably cheaper.
%
%

To maintain concreteness throughout our development, 
we focus on parametric imaging specifically for diffuse optical tomography (DOT);
however, the framework we develop applies to many other nonlinear inversion
problems, such as electrical impedance or resistance tomography (EIT/ERT).
 We observe first that function and Jacobian evaluations in this inversion problem correspond
to evaluations at selected complex frequencies, of a system-theoretic
\emph{frequency response function} of a
parametrized dynamical system,
together with its gradient.
This observation immediately motivates the use of interpolatory parametrized
reduced models as inexpensive surrogates for the full order forward model. For the DOT problem as expressed by the optimization problem (\ref{eq:GenInvProbl}),
these surrogate models are able to approximate both the cost functional and the
associated Jacobian with little loss of accuracy, while significantly reducing
the cost of the overall inversion process.
For the use of model reduction in other optimization applications, we refer
the reader to
\cite{Arian2002,hinze2005proper,Kunisch2008,Antil2011,Antil2012,Druskin2011solution,borcea2012model,
yue2013} and the references therein.

\vspace{0.5cm}
\noindent
In \S\ref{sec:DOTPaLS}, we briefly discuss the medical application
that motivates our work,
and the central role played by
Parametric Level Set (PaLS) parametrization of the medium.
Notably, the PaLS parametrization in effect regularizes the inversion problem, so that no
further regularization is required.
Interpolatory model reduction, as used in the service of
solving the forward problem of DOT-PaLS, is discussed in \S\ref{sec:IntModRed}.
We provide some discussion describing
why our approach appears to work so well for the problems of interest here, and  give
an overview of implementation issues and computational cost.
In \S\ref{sec:NumExp}, we provide four numerical experiments
that demonstrate the effectiveness and accuracy of our
approach.
We offer some discussion of  
future work together with our
concluding remarks  
in \S\ref{sec:Conc}.

\vspace{0.5cm}

\section{Background}\label{sec:DOTPaLS}

\subsection{DOT}
Image reconstruction using DOT typifies the inverse problem
we wish to solve, hence we describe the problem in some detail here.
There are other nonlinear inverse problems with similar structure as well,
such as EIT and ERT.

We assume that the region to be imaged is a rectangular prism (slab):
$\Omega = [-\mathsf{a}_1, \mathsf{a}_1]
\times[-\mathsf{a}_2, \mathsf{a}_2]\times[-\mathsf{a}_3,\mathsf{a}_3]$.
Throughout, $\bfx=(x_1,x_2,x_3)^T$ will refer to spatial location within $\Omega$.
 The top surface of the slab ($x_3 = \mathsf{a}_3$) and the bottom surface ($x_3 = -\mathsf{a}_3$)
 will be denoted as $\partial\Omega_{+}$ and $\partial\Omega_{-}$, respectively.
 The lateral surfaces where either $x_1 = \pm\mathsf{a}_1$ or $x_2 = \pm\mathsf{a}_2$ will
be denoted by $\Gamma$.
Following Arridge \cite{arridge}, we adopt a diffusion model for the photon
flux/fluence $\eta(\bfx,t)$ driven by an input source $g(\bfx,t)$ that is selected
out of a set of $n_{src}$ possible sources, each of which are
(physically) stationary, independently driven, and positioned
on the top surface, $\partial\Omega_{+}$.  We assume there to be functions,
$b_j(\bfx)$, $j=1,\,\ldots,\, n_{src}$ describing the
transmittance field (``footprint") of the $j$th source, so that $g(\bfx,t)=b_j(\bfx)u_j(t)$
for some $j$ and given pulse profile $u_j(t)$.
We assume that observations, $m_i(t)$, are made with a limited number of detectors, say $n_{det}$,
that are presumed to be stationary as well, and located on both the top and bottom surfaces, $\partial\Omega_{\pm}$.   The response characteristics of the sensors are presumed to be captured by functions, $c_i(\bfx)$,  so that $m_i(t)=  \int_{\partial\Omega}c_i(\bfx)\eta(\bfx,t)\,d\bfx$, $i=1,\,\ldots,\, n_{det}$.

The model for the illumination of the region to be imaged then appears as
\begin{align}
  \frac{1}{\nu}\frac{\partial}{\partial t}\eta(\bfx,t) & =
    \nabla \cdot \left(\, D(\bfx) \nabla \eta(\bfx,t) \,\right)
    - \mu(\bfx) \eta(\bfx,t) + b_j(\bfx)u_j(t), \quad \mbox{for }  \bfx\in \Omega ,
\label{eq:ModelPDE1} \\
    0 & = \eta(\bfx,t) +  2\,{\cal A}\,D(\bfx)\,
      \frac{\partial}{\partial \xi}\eta(\bfx,t) , \quad \mbox{for }
      \bfx\in \partial\Omega_{\pm} ,
\label{eq:ModelPDE2} \\
    0 & = \eta(\bfx,t), \quad \mbox{for } \bfx\in \Gamma.
\label{eq:ModelPDE3} \\
    m_i(t) & =  \int_{\partial\Omega} c_i(\bfx)\eta(\bfx,t)\,d\bfx\quad
      \mbox{ for }i=1,\,\ldots,\, n_{det}
\label{eq:ModelPDE4}
\end{align}
(see \cite[p. R56]{arridge}).  $D(\bfx)$ and $\mu(\bfx)$ denote diffusion and absorption
coefficients, respectively;  ${\cal A}$ is
a constant related to diffusive boundary reflection (see \cite[p. R50]{arridge});
$\xi$ denotes the outward unit normal; and $\nu$ is the speed of
light in the medium.

At best, the scalar fields defined by $D(\bfx)$ and $\mu(\bfx)$ are only partially
known. We wish to utilize observations, $\cbfm(t)=(m_1(t),\,m_2(t),\,\ldots,m_{n_{det}}(t))^T$,
 made when the system is
illuminated by a variety of source signals, $\bfu(t)=(u_1(t),\,u_2(t),\,\ldots,u_{n_{src}}(t))^T$, in order
to more accurately determine $D(\bfx)$ and $\mu(\bfx)$.  Accurate determination of
$D(\bfx)$ and $\mu(\bfx)$ is what constitutes ``image reconstruction" for our purposes.
For simplicity,
we assume in our discussion that the diffusivity $D(\bfx)$ is
well specified and that only the absorption field, $\mu(\bfx)$, must be determined.
We also assume that the absorption field, $\mu(\cdot)$, although unknown, is expressible
in terms of a finite set of parameters, $\bfsfp=[p_1,\ldots,\,p_{\ell}]^T$.
An effective parametrization of $\mu(\cdot)$ is fundamental to our undertaking,
and we elaborate on the dependence $\mu(\cdot)=\mu(\cdot,\bfsfp)$
further in \S\ref{ssec:pals}.
In our discussion of the inverse problem and in numerical experiments, we
restrict this model to two dimensions: we consider $\mathsf{a}_2\rightarrow 0$, so $\Omega$ becomes a rectangle in the $x_2$-plane with Dirchlet conditions at
$x_1 = \pm\mathsf{a}_1$ and Robin conditions on the top and
bottom edges where the sources and detectors are located, $x_3 = \pm\mathsf{a}_3$.

A variety of spatial discretizations may be applied to (\ref{eq:ModelPDE1})-(\ref{eq:ModelPDE4}),
finite element methods and finite difference
methods among them, that yield a differential algebraic system of equations represented here as
\begin{equation}  \label{processModelDynSys}
  \frac{1}{\nu}\bfE\,  \dot{\bfy}(t;\bfsfp)  =-\bfA(\bfsfp)\bfy(t;\bfsfp) +
    \bfB\bfu(t)\quad\mbox{with}\quad
    \cbfm(t;\bfsfp) =\bfC\bfy(t;\bfsfp)
\end{equation}
where $\bfy$ denotes the discretized photon flux,
$\cbfm=[m_1,\,\ldots,\, m_{n_{det}}]^T$ is the vector of detector outputs, $\bfC\bfy$ constitutes a
set of quadrature rules for (\ref{eq:ModelPDE4}) applied to the discretized
photon flux; the columns of $\bfB$ are discretizations of the
source ``footprints" $b_j(\bfx)$ for $j=1,\,\ldots,\, n_{src}$;
$\bfA(\bfsfp)=\bfA_{0}+\bfA_{1}(\bfsfp)$ with
$\bfA_{0}$ and $\bfA_{1}(\bfsfp)$ discretizations of the diffusion and absorption
terms, respectively ($\bfA_{1}(\bfsfp)$ inherits the
absorption field parametrization, $\mu(\cdot,\bfsfp)$).
$\bfE$ is generally singular, reflecting the inclusion of the discretized Robin condition (\ref{eq:ModelPDE2}) as an algebraic constraint.

Suppose $\bfY(\omega; \bfsfp)$, $\bfU(\omega)$, and $\bfM(\omega; \bfsfp)$ denote
the Fourier transforms of $\bfy(t;\bfsfp)$, $\bfu(t)$, and $\cbfm(t;\bfsfp)$, respectively.
Taking the Fourier transform of (\ref{processModelDynSys}) and rearranging, we find directly
\begin{equation} \label{FullOrdTransFnc}
  \bfM(\omega;\bfsfp) =\bfPsi\!\left(\omega;\bfsfp\right)\,\bfU(\omega)\quad \mbox{where}\quad
  \bfPsi(\omega;\bfsfp)=\bfC\left(\frac{\imath\;\!\omega}{\nu}\,\bfE\, +\bfA(\bfsfp)\right)^{-1}\bfB,
\end{equation}
where $\omega \in \IR$, and $\bfPsi(\omega;\bfsfp)$ is a mapping from sources (inputs) to measurements (outputs)
in the frequency domain; this is known as the \emph{frequency response function}\footnote{In describing linear dynamical systems, usually the {\it transfer function}
$\bfPsi(s;\bfsfp)=\bfC\left(\frac{s}{\nu}\,\bfE\, +\bfA(\bfsfp)\right)^{-1}\bfB$ is used where  $s\in\IC$ and is not restricted to the imaginary axis. However, for our application here, the measurements are made only on the imaginary axis and it is enough to take $s = \imath \omega$ with $\omega \in \IR$.}
of the dynamical system defined in (\ref{processModelDynSys}).

For any absorption field, $\mu(\cdot,\bfsfp)$, associated with $\bfsfp$,  the vector of
(estimated) observations for the
$i$th input source at frequency $\omega_j$, as predicted by the forward model in
the frequency domain, will be
denoted as $ \bfM_i(\omega_j;\bfsfp)\in\mathbb{C}^{n_{det}}$.
If we stack all the predicted observation vectors for the $n_{src}$ sources and $n_{\omega}$ frequencies, we obtain:
%
\[
\bbM(\bfsfp) = [  \bfM_1(\omega_1;\bfsfp)^T,\, \ldots,\,\bfM_1(\omega_{n_\omega};\bfsfp)^T,\,\bfM_2(\omega_1;\bfsfp)^T, \ldots,
\bfM_{n_{src}}(\omega_{n_\omega};\bfsfp)^T ]^T,~~~
\]
%
which is a (complex) vector of dimension $n_{det}\cdot n_{src} \cdot n_{\omega} $.    We construct the corresponding empirical data vector, $\mathbb{D}$,
from acquired data.
The main optimization problem that must be solved is (cf. (\ref{eq:GenInvProbl})):
\[ \min_{\bfsfp \in \mathbb{R}^\ell} \| \bbM(\bfsfp) - \mathbb{D} \|_2 \]

\subsection{PaLS} \label{ssec:pals}

The concept of level sets was first introduced by Osher and
Sethian in \cite{osher1988fronts} and has since gained momentum; see \cite{santosa1996level,dorn2006level,burger2005survey,DoelAsch06,DoelAsch07,DoelAsch10}.
Traditional level set methods need specialized optimization and
well-honed regularization to overcome frequent sensitivity
to ill-posedness or noise in the problem.

\begin{figure}[t]
\centering
\begin{tabular}{cc}
\includegraphics[width=2.5in,height=2in]{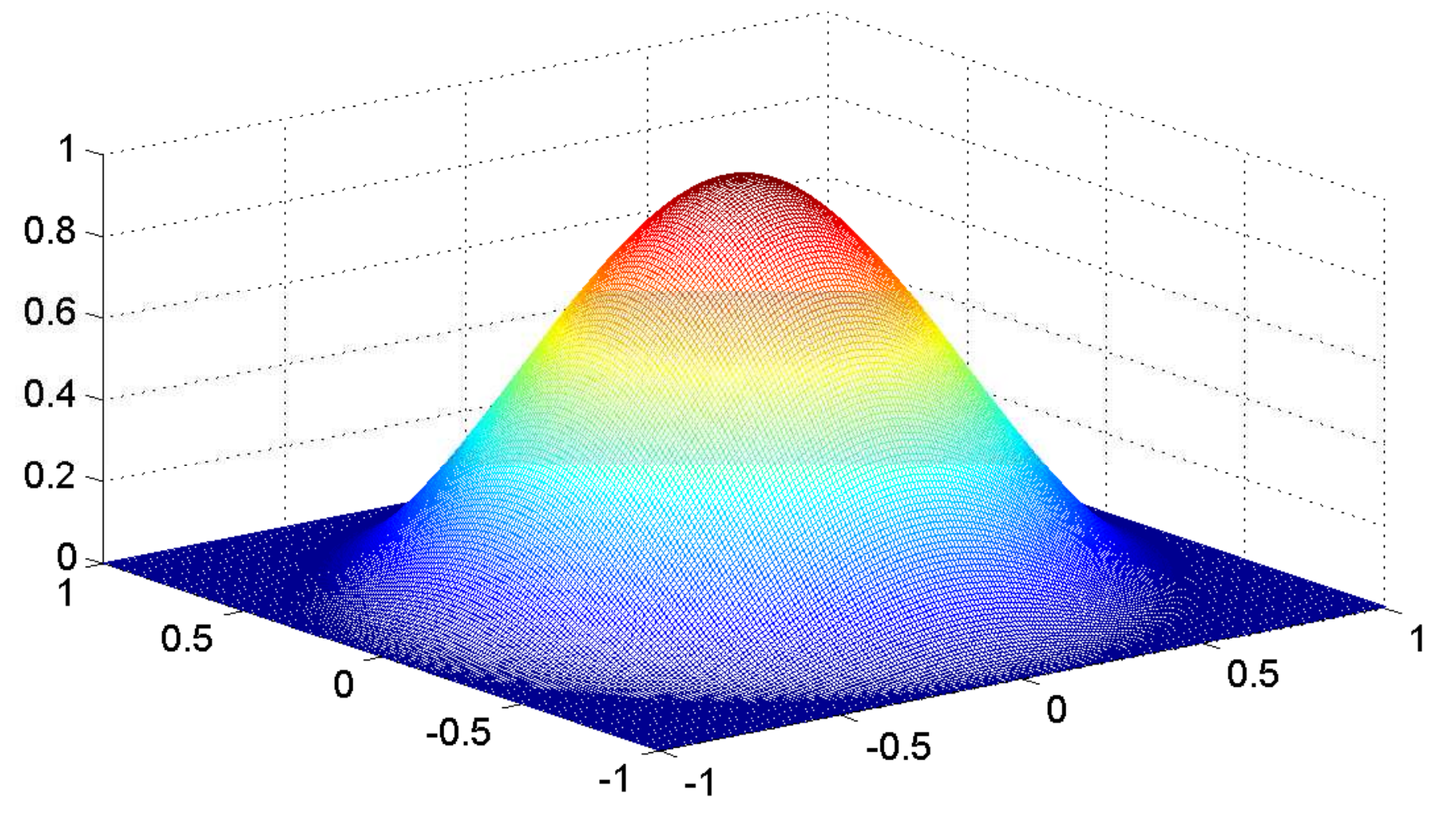} &
\includegraphics[width=2.5in,height=2in]{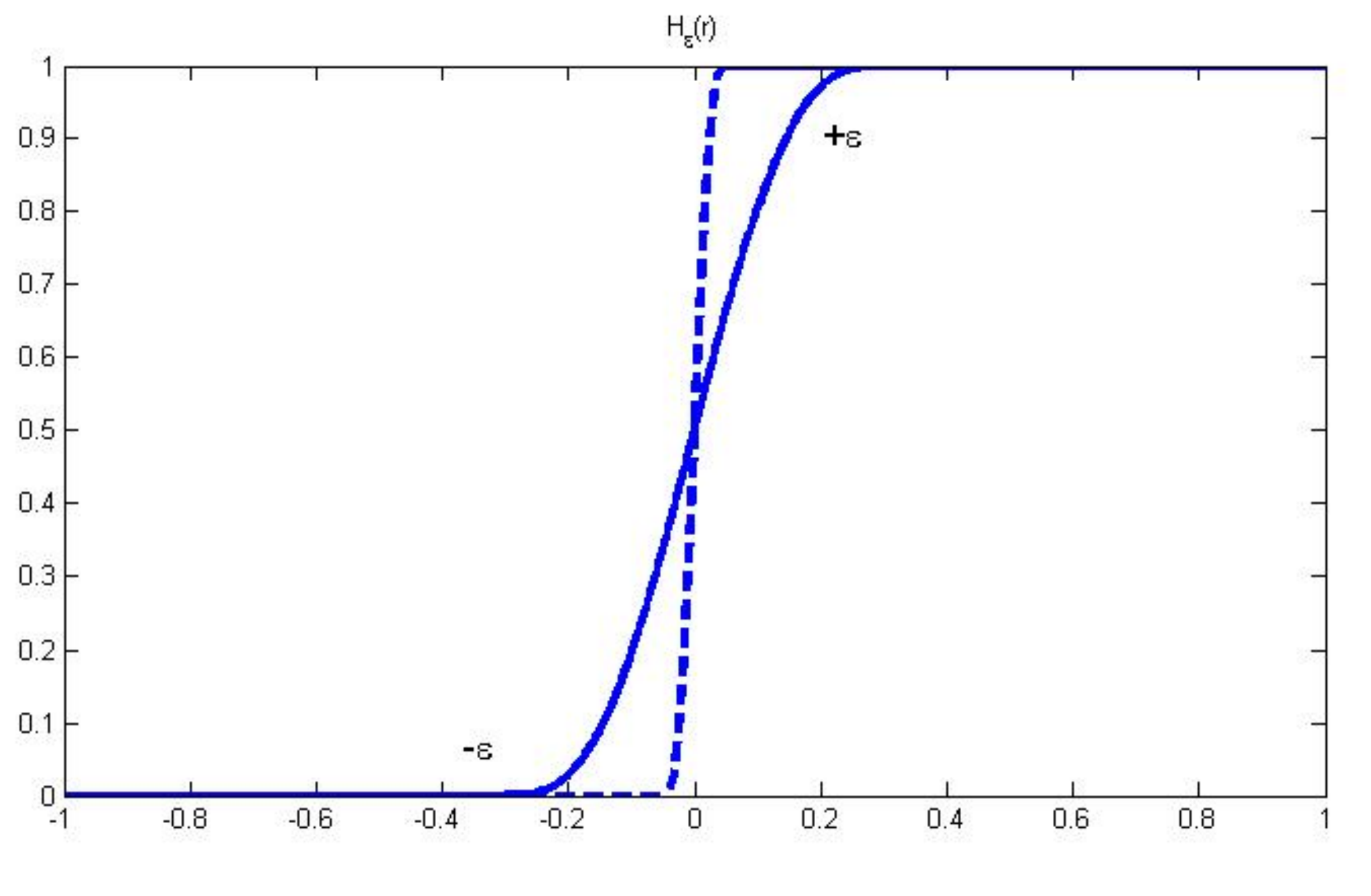} \\
(a) & (b)
\end{tabular}
\caption{\label{fig:csrbf} (a) The CSRBF is given by
$\varphi(r) = \left( \max(0,1-r) \right)^2(2r + 1)$ with $r = \sqrt{x^2 + y^2}$ \cite{wendland2005scattered}.
(b) Two graphical representations of $H_{\epsilon}(r)$,
our approximation~to~a~Heaviside~function.  The width of the transition
region is defined by $\epsilon$.  The dashed line shows an $H_{\epsilon}(r)$ for
a larger value of $\epsilon$.}
\end{figure}

The use of parametric level set (PaLS) representations of images was recently shown to be
beneficial in a large number of inverse problems \cite{Aghasi_etal11}, in terms of
drastically reducing the dimension of the search space (over typical
voxel-based formulations) and providing implicit regularization for
a suitable number of basis functions in the representation
(as opposed to Tikhonov regularization where one must handle the difficult
problem of choosing the regularization parameter).
We use compactly supported radial basis functions (CSRBF)
for the PaLS representation, as in \cite{Aghasi_etal11}.
An example of a CSRBF is given in Figure \ref{fig:csrbf}.
See {\cite{wendland2005scattered} for other choices.

Let $\varphi: \mathbb{R}^{+} \rightarrow \mathbb{R}$ denote a CSRBF;
define $ \| \bfx \|^{\dagger} := \sqrt{ \|\bfx \|_2^2 + \gamma^2 }$; and
\[
  \phi(\bfx,\bfsfp_{-}) := \sum_{j=1}^{m_0} \alpha_j
  \varphi( \| \beta_j( \bfx - \bfchi_j) \|^{\dagger} ),
\]
where $\bfsfp_{-}$ denotes a vector of unknown parameter values comprised
of the expansion coefficients $\alpha_j$, dilation factors $\beta_j$,
and center locations $\bfchi_j$.  The CSRBF $\varphi$ is assumed to be
sufficiently smooth
and the value of $\gamma$ is small and positive; $\gamma$  is introduced here
to avoid complications of derivatives of the CSRBF parameters at zero.
We assume  the absorption field is compatible with a level-set parametrization
\begin{equation} \label{eq:ff}
  \mu(\bfx,\bfsfp) = \mu_{in}(\bfx) H_{\epsilon}(\phi(\bfx,\bfsfp_{-}) - c)
    + \mu_{out}(\bfx) (1 - H_{\epsilon}(\phi(\bfx,\bfsfp_{-}) - c) ),
\end{equation}
where $H_{\epsilon}(r)$ denotes a continuous approximation to the Heaviside
function, shown in Figure~\ref{fig:csrbf}, $c$ is the height of the level set of interest,
and $\bfsfp$ represents the concatenation of $\bfsfp_{-}$ and any additional
parameters that define $\mu_{in}(\bfx)$ and $\mu_{out}(\bfx)$.
The effect is that $\mu(\bfx,\bfsfp)$ is very nearly piecewise-constant,
having the value $\mu_{in}(\bfx)$ if $\bfx$ is inside the region defined
by the c-level set, and having value $\mu_{out}(\bfx)$ otherwise.
In all that follows, we assume $\bfsfp=\bfsfp_{-}$ for ease of discussion,
since incorporating parameters for $\mu_{in}$ and $\mu_{out}$ into the
reduced order model (ROM) framework is straightforward.
We are able to model nonconstant diffusion similarly.
Here, we assume that diffusion is known for ease of exposition.

This PaLS formulation can capture edges and complex boundaries
with relatively few basis functions \cite{Aghasi_etal11}.
The compact support of the basis functions is an important
advantage for nonlinear optimization,
since not all parameters may need to be updated at each iteration (see \cite{Aghasi_etal11}).
The TREGS (pronounced \emph{t\={e}~reks}) 
method \cite{StuKil11c} has proved to
be fast and reliable at solving the nonlinear least squares problem
for the parameter vector describing the
absorption images, and we therefore use this algorithm for all our
numerical results.  The interested reader is directed to \cite{StuKil11c}
for further details of the optimization algorithm.

\section{Interpolatory Model Reduction for DOT-PaLS}\label{sec:IntModRed}

\subsection{A systems theoretic perspective}
In the course of finding a minimizer to (\ref{eq:GenInvProbl}), a single evaluation of $\bbM(\bfsfp) - \mathbb{D}$
involves computing (for all $i$ and $j$)
\eqs \label{eq:CYi_om_j}
  \bfM_i(\omega_j; \bfsfp)= \bfC\bfY_i(\o_j;\bfsfp) =
    \bfC\left(\imath\frac{\omega_j}{\nu}\,
    \bfE\, +\bfA(\bfsfp)\right)^{-1}\bfB\,\bfU_i(\omega_j)
    = \bfPsi( \omega_j) \bfU_i(\omega_j) , \quad
\eqe
where $\bfPsi(\omega;\bfsfp)$ is the frequency response function defined in (\ref{FullOrdTransFnc})
and $ \bfU_i(\omega_j) = \mathbf{e}_i$,  the $i$th
column of the identity matrix; $\bfU_i(\omega_j)$ excites the $i$th source location with
a pure sinusoid signal at frequency $\omega=\omega_j$.
From a systems theoretic perspective,  the objective function evaluations
are precisely evaluations of a frequency response function on the imaginary
axis. These function evaluations require solving
large, sparse linear systems
$\left(\imath\frac{\omega_j}{\nu}\,\bfE\, +\bfA(\bfsfp)\right)
\bfv = \bfB\,\mathbf{e}_i$, for all source locations $\mathbf{e}_i$
and frequencies $\omega_j$.
One major advantage of the PaLS parametrization is that the work for
a function evaluation is  independent of the number of parameters
used, although it does depend on the order of the system and the numbers of sources and detectors, so the computational complexity will still be quite high.  Consider,
for example, an $8\mbox{cm} \times 8\mbox{cm} \times 4 \mbox{cm}$ region
with detectors spaced $0.5$mm apart and sources spaced 2mm apart on the top and
bottom (consistent with current systems).
{\em One function evaluation requires $40 \times 40 \times 2$ linear systems
with about $2\times {10^6}$ unknowns for each frequency}.

Similar observations hold for the Jacobian computation as well.
The Jacobian is constructed using an adjoint-type (or co-state)
approach that exploits the fact that the number of detectors
is roughly equal to the number of sources,
as discussed in \cite{HabAsch00} and \cite[p. 88]{Vogel}.
Using (\ref{FullOrdTransFnc}) and (\ref{eq:CYi_om_j})
and differentiating $ \bfM_i(\omega_j;\bfsfp)$
with respect to $\bfsfp$ gives
\begin{equation} \label{eq:jaccalc}
  \frac{\partial }{\partial p_k} \bfM_i(\omega_j;\bfsfp)= \frac{\partial }{\partial
    p_k}\left[\bfPsi\!\left(\omega_j;\bfsfp\right)\right]\,\bfU_i(\omega_j)
    =  - \bfZ(\omega_j;\bfsfp)^T\,
    \frac{\partial }{\partial p_k}\bfA(\bfsfp)\, \bfY_i(\omega_j; \bfsfp),
\end{equation}
where, for each $\omega_j$ and any $\bfsfp$ we can compute $\bfZ(\omega_j;\bfsfp)$ from
$ \left(\imath\frac{\omega_j}{\nu}\,\bfE\, +\bfA(\bfsfp)\right)^T \bfZ(\omega_j;\bfsfp) =  \bfC^T .$
 Equation (\ref{eq:jaccalc}) reveals that
Jacobian evaluations in this problem correspond to evaluation of
partial derivatives of the transfer function $\bfPsi(s,\bfsfp)$ with respect to the parameters.
The matrices $\frac{\partial }{\partial p_k}\bfA(\bfsfp)$ need to be computed
only once.
Thus, the computational cost for evaluating
the Jacobian with respect to $\bfsfp$, apart from computing $\bfM_i(\omega_j; \bfsfp)$ (which is also necessary for the function evaluation),
consists mainly of the cost of computing $\bfZ(\omega_j;\bfsfp)$, that is, solving $n_{det}\cdot n_{\omega}$ linear systems of order equal to the
number of degrees of freedom.

In short, the solution of the nonlinear least squares problem (\ref{eq:GenInvProbl})
requires evaluating both $\bfPsi(\omega;\bfsfp)$ and $\nabla_{\bfsfp}\bfPsi(\omega;\bfsfp)$
for many values of $\bfsfp$ and $\omega$ and hence solving
a substantial number of large, sparse linear systems.
This is a critical bottleneck, and
effective strategies pivot on efficient solution of the forward problem.
Our remedy for this bottleneck has many features in common with methods
of \emph{model reduction} (see \cite{antoulas2005approximation,Ant2010imr}):
we seek a frequency response function $\widehat{\bfPsi}(\omega;\bfsfp)$
that is easy to evaluate, yet provides a high-fidelity approximation to
$\bfPsi(\omega;\bfsfp)$  over parameters and  frequencies of interest.
Likewise, we require that $\nabla_{\bfsfp}\wh{\bfPsi}(\omega;\bfsfp)$ is easy
to evaluate and that
$\nabla_{\bfsfp}\bfPsi(\omega;\bfsfp)\approx \nabla_{\bfsfp}\widehat{\bfPsi}(\omega;\bfsfp)$
over the same range of arguments.

\subsection{Surrogate Forward Model via Parametric Model Reduction}
Recall  the (time-domain) dynamical systems representation
(\ref{processModelDynSys}) of the
``original" parametric system, i.e., the forward model. With this in mind,
we are seeking a much smaller parametric model, of order $r \ll n$,
that is able to replicate the input-output map
of the original model (\ref{processModelDynSys}).
So, we construct the small dynamical system
\vspace{-1ex}
\begin{equation} \label{eq:ParamSysROM}
 \frac{1}{\nu}\,\widehat{\bfE}\, \dot{\widehat{\bfy}}_r(t;\bfsfp)  =
  - \widehat{\bfA}(\bfsfp)\, \widehat{\bfy}(t;\bfsfp) + \widehat{\bfB}\, \bfu(t)   \quad\mbox{with}\quad  \widehat{\boldsymbol{m}}(t;\bfsfp)  =
   \widehat{\bfC}\, \widehat{\bfy}(t;\bfsfp) ,
\end{equation}
where the new state vector
$\widehat{\bfy}(t;\bfsfp)\in \IR^r$,
$\widehat{\bfE},\widehat{\bfA}(\bfsfp) \in \IR^{r \times r}$,
$\widehat{\bfB} \in \IR^{r \times n_{src}}$, and  $ \widehat{\bfC}\in \IR^{n_{det}\times r}$
such that ${\boldsymbol{m}}(t;\bfsfp) \approx \widehat{\boldsymbol{m}}(t;\bfsfp)$.
The surrogate (reduced) frequency response function appears as
\eqs \label{eq:Psi_hat_1}
  \widehat{\bfPsi}(\omega;\bfsfp) & = &
    \widehat{\bfC} \left(\frac{\imath\,\!\omega}{\nu}  \widehat{\bfE} + \widehat{\bfA}(\bfsfp)\right)^{-1} \widehat{\bfB}.
\eqe
Since $\widehat{\bfB}$ has the same number of columns as $\bfB$
and $\widehat{\bfC}$ has the same number rows as $\bfC$, $\bfPsi(\omega;\bfsfp)$ and
$\widehat{\bfPsi}(\omega;\bfsfp)$ have the same row and column dimension,
although the state vector $\widehat{\bfy}(t;\bfsfp)$  occupies a much
lower dimensional space than does the original $\bfy(t;\bfsfp)$.
However, evaluating $\bfPsi(\omega;\bfsfp)$ requires the solution of linear systems
of dimension $n$, whereas evaluating $\widehat{\bfPsi}(\omega;\bfsfp)$ requires
the solution of linear systems only of dimension $r \ll n$. This will reduce
the cost of the forward problem drastically. A similar discussion applies to evaluating
$\nabla_{\bfsfp}\bfPsi(\omega;\bfsfp)$ vs $\nabla_{\bfsfp}\widehat{\bfPsi}(\omega;\bfsfp)$.

We construct the surrogate parametric model using projection.
Suppose full rank matrices $\bfV \in \IC^{n \times r}$ and $\bfW\in \IC^{n \times r}$
are specified. Using the approximation ansatz
that the full state $\bfy(t;\bfsfp)$ evolves roughly near the $r$-dimensional
subspace $\Ra{\bfV}$, we write
$\bfy(t;\bfsfp) \approx \bfV \widehat{\bfy}(t;\bfsfp)$ and enforce the Petrov-Galerkin condition
$$
  \mathbf{W}^{T}\left(\frac{1}{\n}\bfE \bfV\dot{\widehat{\bfy}} (t) +
    \bfA(\bfsfp)\bfV\widehat{\bfy}(t)-\bfB\,\bfu(t)\right)=\mathbf{0}, \qquad
    \quad \widehat{\boldsymbol{m}}(t) = \bfC\bfV \widehat{\bfy}(t).
$$
\noindent
to obtain the reduced-system given in (\ref{eq:ParamSysROM}) with
the reduced matrices given by
\begin{align} \label{eq:PGProj}
  \widehat{\bfE} = \bfW^{T}\bfE\bfV, \quad
    \widehat{\bfA}(\bfsfp) = \bfW^{T} \bfA(\bfsfp) \bfV,
    \quad \widehat{\bfB}= \bfW^{T}\bfB, \quad \widehat{\bfC} = \bfC\bfV.
\end{align}
The computational advantages of using the surrogate model $\widehat{\bfPsi}(\omega;\bfsfp)$
in place of the full-model ${\bfPsi}(\omega;\bfsfp)$ will be explained
in more detail in \S\ref{sec:IntPMOR}, after we explain how we choose
model reduction bases.

\subsection{Interpolatory Parametric Model Reduction}  \label{sec:IntPMOR}
Several parametric model reduction methods exist to select
$\bfV \in \IC^{n \times r}$ and $\bfW \in \IC^{n \times r}$; see, e.g., \cite{prud2002reliable,RozHP08,BuiThanh2008,nguyen2008best,BauB09,haasdonk2011erm,Hay09,veroy2003posteriori,BuiThanh2008_AIAA}
and the references therein. Recall that our function and Jacobian evaluations
correspond to transfer function evaluations and their derivatives.
Thus, given a parameter vector $\hat{\bfsfp} \in \IR^\ell$ and a frequency $\hat{\omega}\in \IR$,
we would like to use a reduced parametric model of the form
(\ref{eq:ParamSysROM})
whose frequency response function $\widehat{\bfPsi}(\omega;\bfsfp)$ satisfies
%
\begin{equation} \label{eq:hermite}
  \widehat{\bfPsi}\left( \hat{\omega};\hat{\bfsfp}\right) =
    \bfPsi\left( \hat{\omega};\hat{\bfsfp}\right)
    ~~~\mbox{and}~~~
    \nabla_{\bfsfp}\widehat{\bfPsi}\left( \hat{\omega};\hat{\bfsfp}\right) = \nabla_{\bfsfp}\bfPsi\left(\hat{\omega};\hat{\bfsfp}\right).
\end{equation}
This perfectly fits in the context of interpolatory parametric model reduction
(\cite{bond2005pmo,BauBBG09,Daniel2004,gunupudi2003ppt,FenB08}).
The following theorem, from \cite{BauBBG09,Ant2010imr}, shows how to construct
the reduction bases $\bfV$ and $\bfW$, so that the reduced transfer
function satisfies (\ref{eq:hermite}). In \cite{BauBBG09,Ant2010imr},
all the matrices in the original system are assumed to have parametric
dependency, not only $\bfA(\bfsfp)$. We present the theorem in
the context of the DOT-PaLS problem.
\begin{theorem} \label{thm:parammor}
Suppose $\bfA(\bfsfp)$ is continuously differentiable in a neighborhood
of $\hat{\bfsfp} \in \IR^{\ell}$. Let $\hat{\omega} \in \IC$,
and both
$\frac{\imath\,\!\hat{\omega}}{\nu}\,\bfE +\bfA(\hat{\bfsfp})$ and
$\frac{\imath\,\!\hat{\omega}}{\nu}\,\widehat{\bfE} +\widehat{\bfA}(\hat{\bfsfp})$ be invertible.
$$  \mbox{\textit{If}}
	\left(\frac{\imath\,\!\hat{\omega}}{\nu}\,\bfE
     +\bfA(\hat{\bfsfp})\right)^{-1}\bfB\in
    \Ra{\bfV}~\mbox{ \textit{and} }~
	 \left(\bfC\left(\frac{\imath\,\!\hat{\omega}}{\nu}\,\bfE
     +\bfA(\hat{\bfsfp})\right)^{-1}\right)^{T} \in \Ra{\bfW},
$$
then, the reduced parametric model of (\ref{eq:PGProj}) satisfies
$$
	  \bfPsi\!\left(\hat{\omega},\hat{\bfsfp}\right) =\widehat{\bfPsi}\!\left(\hat{\omega},\hat{\bfsfp}\right),~~
    \nabla_{\bfsfp}\bfPsi\!\left(\hat{\omega},\hat{\bfsfp}\right)  =\nabla_{\bfsfp}\widehat{\bfPsi}\!\left(\hat{\omega},\hat{\bfsfp}\right), ~~ \mbox{\textit{and}} ~~ \bfPsi'\!\left(\hat{\omega},\hat{\bfsfp}\right) =\widehat{\bfPsi}'\!\left(\hat{\omega},\hat{\bfsfp}\right),
$$
where $'$ denotes the derivative with respect to $\omega$.
\end{theorem}

So, we can match both function and gradient values exactly without explicitly computing them,
requiring only that $\Ra{\bfV}$ and $\Ra{\bfW}$ contain particular vectors,
see \cite{BauBBG09,Ant2010imr}.  Thus, by constructing
$\bfV$ and $\bfW$ as in
 Theorem \ref{thm:parammor}, we satisfy the
desired interpolation conditions (\ref{eq:hermite}) for the inverse
problem.
This means that for the frequency and parameter interpolation
points selected to construct $\bfV$ and $\bfW$, an optimization algorithm
based on the Gauss-Newton approach (such as TREGS)
would proceed identically for the reduced forward model as for the full-order forward
model. In other words,
if, by chance, the optimization algorithm were to pick only parameter vectors
that are in the set of interpolations points used to construct $\bfV$ and $\bfW$,
there would be no difference between using the full forward model and the
surrogate forward model.

We  are able to avoid recomputing the projection bases, $\bfV$ and
$\bfW$, as the optimization process generates new parameter vectors.
Instead, we use parametrized reduced
order models that can be updated (relatively) cheaply for each new parameter
vector and the same projection bases are reused (typically) in each step. (see \S\ref{ssec:GlobalBases}).
Projection bases are constructed that are able to provide
parametrized reduced models that retain high fidelity over the range
of parameter values that are needed to capture features
of the full forward model that emerge in the course of optimization.
This approach ensures that optimization using the
\emph{surrogate} model will provide solutions of the same quality
as would occur had the full model been used instead.
As will be discussed below and demonstrated in the numerical experiments,
this turns out to be possible at very low cost.


\subsubsection{Using Global Basis Matrices for Projection} \label{ssec:GlobalBases}
In this paper, we use a \emph{global basis} approach to construct the model
reduction bases: we construct two constant
matrices $\bfV$ and  $\bfW$, built by sampling at multiple parameter values, that
capture sufficient information to produce near interpolants across the
needed range of parameter values.
This approach contrasts with \emph{local bases methods} where $\bfV$ and a $\bfW$
must be updated as $\bfsfp$ varies;
see, e.g., \cite{Amsallem2008,Panzer_etal2010,AmsallemFarhat2011,Degroote2010},
and the references therein.

Given parameter sample points
$\boldsymbol{\pi}_1,\ldots,\boldsymbol{\pi}_K$,
following Theorem \ref{thm:parammor}, we construct, for $i=1,\ldots,K$,
the local basis matrices
\eqs   \label{eq:localVi}  
  \bfV_i & = & \left[ \bfV_{i,1}, \bfV_{i,2}, \ldots, \bfV_{i,n_\o}\right]   \\
  & = &
  \left[ \left(\imath \frac{\omega_1}{\nu}\bfE + \bfA(\boldsymbol{\pi}_i)\right)^{-1}\bfB,\ldots, \left(\imath \frac{\omega_{n_\o}}{\nu}\bfE + \bfA(\boldsymbol{\pi}_i)\right)^{-1}\bfB \right] ,  \nonumber
\eqe
where
$\bfV_{i,j} = \left(\imath \frac{\o_j}{\n}\bfE + \bfA(\boldsymbol{\pi}_i)\right)^{-1}\bfB$
and
$\omega_1,\ldots,\omega_{n_\omega}$ are the frequency interpolation
points. Note that the frequency interpolation points in this
context are predetermined by the experimental set-up, and we do
not need to search for optimal frequency interpolation, as usually
required in rational interpolation based model reduction;
see, e.g., \cite{gugercin2008hmr,BauBBG09}.
Since the local basis matrices might have common components among each other, we
eliminate those components by taking the left singular vectors of the concatenated local basis
matrices $ [\bfV_1, \bfV_2,\ldots, \bfV_{K}]$ as
the global basis $\bfV$; thus the reduced order $r$ becomes the number of non-zero singular values.
Similar steps are applied to construct $\bfW$,
but using the adjoint system with columns of  $\bfC^T$ on the right-hand side
(see Theorem \ref{thm:parammor}).
We are then assured that the reduced parametrized model
with $\widehat{\bfE}= \bfW^{T}\bfE\bfV$,
$\widehat{\bfA}(\bfsfp) = \bfW^{T} \bfA(\bfsfp) \bfV$,
$\widehat{\bfB}= \bfW^{T}\bfB$, and $\widehat{\bfC} = \bfC\bfV$
will satisfy
(\ref{eq:hermite}) at {\em every}
$(\omega,\bfsfp) =  (\omega_j,\boldsymbol{\pi}_i)$ for
$j=1,\ldots,n_\omega$ and $i=1,\ldots,K$.

Now, consider the efficient evaluation of $\widehat{\bfPsi}(\omega;\bfsfp) =
\widehat{\bfC} \left(\frac{\imath\,\!\omega}{\nu}  \widehat{\bfE} +
\widehat{\bfA}(\bfsfp)\right)^{-1}\widehat{\bfB}$ given
global projection bases $\bfV$ and $\bfW$.
Since $\widehat{\bfE}$, $\widehat{\bfB}$ and $\widehat{\bfC}$ are fixed,
we need to consider only the efficient computation of $\widehat{\bfA}(\bfsfp)$,
the cost of solving small $r \times r$ systems being negligible.
In our current setting,
$\bfA(\bfsfp) = \bfA_0 + \bfA_1(\bfsfp)$,
where $\bfA_0$ is constant and $\bfA_1(\bfsfp)$ carries the parametric dependency.
Hence, $\widehat{\bfA}(\bfsfp) = \bfW^T \bfA_0 \bfV +  \bfW^T\bfA_1(\bfsfp) \bfV$,
and $\bfW^T \bfA_0 \bfV$ can be precomputed.
Only $\bfW^T\bfA_1(\bfsfp) \bfV$ needs to be recomputed after
updating the parameter vector.
Notice that this holds even if the diffusion is also a function of $\bfsfp$.

Assume that a finite difference discretization is used for the underlying PDE.
In that case, $\bfA_1(\bfsfp)$ is diagonal,
making $\bfA_1(\bfsfp) \bfV$ very cheap.
So, the main cost is the computation of
in $\bfW^T (\bfA_1(\bfsfp) \bfV)$.
If diffusion also depends on $\bfsfp$, the multiplication $\bfA_1(\bfsfp) \bfV$
will be $7$ times as expensive as multiplying $\bfV$ by a diagonal
matrix in the $3$D case
and only $5$ times as expensive in the $2$D case.
So, the computation of $\bfW^T\bfA_1(\bfsfp) \bfV$ is very cheap
compared with the solution of many large, sparse linear systems
for the full forward model.
Similar computational costs will occur with finite element discretizations,
although in that case we would not have a diagonal matrix
even when inverting only for the absorption.

Further cost reductions are possible due to the approximate Heaviside
function, $H_\epsilon$, in (\ref{eq:ff}).
When the change in the parameter vector is small, which is typically
enforced step-wise in nonlinear optimization, only a small subset of
the coefficients of $\bfA_1(\bfsfp)$ will change.
We can exploit this in the update of $\bfW^T\bfA_1(\bfsfp) \bfV$.
Let $\bfDe = \bfA_1(\bfsfp_2) - \bfA_1(\bfsfp_1)$; typically
$\bfDe_{k,k} = 0$ for all but a modest number of entries, say,
$\bfDe_{k_1,k_1}, \ldots, \bfDe_{k_q,k_q}$.
We have
$\widehat{\bfA}(\bfsfp_2) = \widehat{\bfA}(\bfsfp_1) + \bfW^T \bfDe \bfV$,
where only $\bfW^T \bfDe \bfV$ must be computed.
This requires the multiplication of
a modest number of entries of $\bfW$ and $\bfDe \bfV$
with computational cost of only
$O(r^2 q)$ flops, a constant cost in terms of $n$.
Again, this would also be true if the diffusion depended
on $\bfsfp$ or when using finite element discretization,
except that $\bfDe$ would not be diagonal (but still mostly zero).

The Jacobian can be evaluated cheaply in a similar way.
For the $k$th parameter $p_k$, we obtain
\eqsn
  \nabla_{p_k}\widehat{\bfPsi}(\omega;\bfsfp)= -\widehat{\bfC}
    \left( \frac{\imath\,\!\omega}{\nu} \widehat{\bfE} + \widehat{\bfA}(\bfsfp) \right)^{-1}
    \left( \frac{\partial}{\partial p_k} \widehat{\bfA}(\bfsfp) \right)
    \left( \frac{\imath\,\!\omega}{\nu}  \widehat{\bfE} + \widehat{\bfA}(\bfsfp) \right)^{-1} \widehat{\bfB} ,
\eqen
where
\[
  \frac{\partial}{\partial p_k}\widehat{\bfA}(\bfsfp) =
    \bfW^T \left( \frac{\partial}{\partial p_k} \bfA_1(\bfsfp)\right)  \bfV ,
\]
which can be computed efficiently due to the structure explained above.


\subsection{Analysis of Global Bases Approach}  \label{sec:anal_gb}

In this section, we motivate why the approach discussed above
is likely to be very effective
for our problem and similar problems,
leaving details and proofs to a follow-up paper.
In particular, we argue that for our problem
the global projection bases for the reduced
model, (\ref{eq:localVi}) for $\bfV$ and similar for $\bfW$,
are not expected to change
much as a function of the parameters.
This means that (1) we need to compute local bases only
for a modest number of interpolation points to obtain a good approximation of the full forward model,
thus keeping both the number of full size linear systems to solve and the
size of the reduced model small to modest;
(2) it might be possible to compute the global projection bases
once (off-line) and use them for a range of distinct image reconstructions;
and (3) the difference in optimization using the full model
and using the reduced model is small, leading to an effective optimization.
In Section~\ref{sec:NumExp}, we demonstrate that at least
for our proof-of-concept problems, this indeed seems to be the case.

For ease of discussion, 
we consider only finite difference discretizations and
invert only for absorption (not an unreasonable
assumption for DOT).   For simplicity, we assume our model uses an exact Heaviside
function in (\ref{eq:ff}).  Typically, $|\omega|$ is small relative to $\n$ and
the mesh width $h$ is small.

%
%
%
%

Under these assumptions, the matrix
$(\imath \,\!\omega/\n) \bfE + \bfA(\bfsfp) = \bfA_0 + (\imath \,\!\omega/ \n)\bfE + \bfA_1(\bfsfp)$ varies
only on the diagonal, and changes are quite small relative to the
magnitude of the matrix coefficients as $\bfsfp$ changes.
The matrix $\bfA_0$ represents the discretization of $-\nabla \cdot(D(x)\nabla \h)$
for a particular choice of diffusion field $D(x)$
scaled by $h^2$, so that the diagonal coefficients
of $\bfA_0$ are $O(1)$.
The matrix $\bfA_1(\bfsfp)$ represents the discretization of the
absorption term $\m(\bfx;\bfsfp)\h$ scaled by $h^2$,
and $\bfE$ is the identity scaled by $h^2$,
except for the diagonal coefficients corresponding to the
top and bottom grid points, which are zero.
Hence, we can write
\eqs \label{eq:-A0+h2D}
  \frac{\imath \,\!\omega}{\n} \bfE + \bfA(\bfp) = \bfA_0 + h^2 \bfD(\omega,\bfp) ,
  \quad \mbox{where} \quad
  \bfD(\omega,\bfp) = \frac{\imath \,\!\omega}{\n}\bfE + \bfA_1(\bfsfp)
\eqe
and the coefficients $\bfD_{jj}$ are $O(1)$.
Suppose $\m_1$ corresponds to the absorption coefficient in normal tissue and $\m_2$
corresponds to the absorption coefficient in anomalous tissue, respectively
$\m_{out}$ and $\m_{in}$ in (\ref{eq:ff}).
Let $\bfI_1(\bfp)$ be a diagonal matrix representing the discrete indicator,
or characteristic function, for normal tissue, meaning it contains zeros on the diagonal
for indices corresponding to pixels with anomalous tissue and ones on the
diagonal otherwise. Then define $\bfI_2(\bfp) = \bfI - \bfI_1(\bfp)$ as the diagonal matrix representing
the indicator function for anomalous tissue.  Note that in general we expect that $\bfI_2(\bfp)$ has mostly
zeros on the diagonal, reflecting the assumption of a relatively small anomaly in mostly healthy tissue.
Using the definitions above, we may
decompose $\bfD(s,\bfp)$ as
\eqs \label{eq:Dsp}
  \bfD(\omega,\bfp) = \frac{\imath \,\!\omega}{\n} \bfE + \m_1 \bfI_1(\bfp) + \m_2 \bfI_2(\bfp) ,
\eqe

Now, we can show that the matrices $\bfV_{i,k}$ in (\ref{eq:localVi}), where
$\bfsfp_k$ is the approximate solution at step $k$ in the nonlinear optimization,
will generally stay close to each other. This suggests that global
projection bases computed, for example, using the parameter vectors from the
first few iterations are close to local projection bases
that could be computed for other parameter vectors later in the optimization.
Hence, the reduced transfer function remains accurate for the duration
of the optimization. We show the close proximity of these spaces
numerically in Figure~\ref{fig:MinSinVal} for two test problems discussed
in more detail in the next section. This also suggests that the
computation of the reduced transfer function is not sensitive to the
choice of interpolation points in parameter space.

First, consider the matrices
$\bfV_{i,1} = (\frac{\imath \,\!\omega_1}{\nu} \bfE + \bfA(\bfsfp_i))^{-1} \bfB$, for $i = 1, \ldots, K$,
where $\omega_1 = 0$,
\[
  \bfV_{i,1} = \bfA(\bfsfp_i)^{-1}\bfB
         = \left( \bfA_0 + h^2 \m_1\bfI + h^2(\m_2 - m_1)\bfI_2(\bfsfp_i) \right)^{-1} \bfB ,
\]
where we have used the fact that
$\m_1 \bfI_1(\bfsfp) + \m_2 \bfI_2(\bfsfp) =
\m_1 \bfI + \m_2 \bfI_2(\bfsfp) - \m_1 \bfI_2(\bfsfp)$.
We attempt to show that $\Ra{\bfV_{i,1}}$ is close to
$\Ra{ ( \bfA_0 + h^2 \m_1 \bfI)^{-1} \bfB }$, where
$\bfA_0 + h^2\m_1 \bfI$ is the operator for constant absorption of healthy tissue.
Factoring $\bfA(\bfsfp_i)$ gives
\begin{align*}
  \left( \bfA_0 + h^2 \m_1\bfI + h^2(\m_2 - \m_1)\bfI_2(\bfsfp_i) \right) = & \\
  \left( \bfA_0 + h^2 \m_1\bfI \right) &
    \left( \bfI + h^2(\m_2 - \m_1) \left( \bfA_0 + h^2 \m_1\bfI \right)^{-1}
    \bfI_2(\bfsfp_i) \right) ,
\end{align*}
and we get
\eqs
\nonumber
  \bfV_{i,1} & = &
  \left[ \left( \bfA_0 + h^2 \m_1\bfI \right)
    \left( \bfI + h^2(\m_2 - \m_1) \left( \bfA_0 + h^2 \m_1\bfI \right)^{-1}
    \bfI_2(\bfsfp_i) \right) \right]^{-1} \bfB \\
\nonumber
  & = &
  \left( \bfI + h^2(\m_2 - \m_1) \left( \bfA_0 + h^2 \m_1\bfI \right)^{-1}
    \bfI_2(\bfsfp_i) \right)^{-1}
    \left( \bfA_0 + h^2 \m_1\bfI \right)^{-1} \bfB \\
\nonumber
  & = &
  \left( \bfI - h^2(\m_2 - \m_1)
    \left[ ( \bfA_0 + h^2 \m_1\bfI )^{-1}\bfI_2(\bfsfp_i) \right] \right. \\
\nonumber
  &&
    \qquad + \,
    h^4 (\m_2 - \m_1)^2 \left[( \bfA_0 + h^2 \m_1\bfI )^{-1}\bfI_2(\bfsfp_i) \right]^2 \\
\label{eq:Vi1}
  &&
  \left.
    \qquad - \, \cdots \right)
    \left( \bfA_0 + h^2 \m_1\bfI \right)^{-1} \bfB .
\eqe
Although the result depends on the rate of convergence of the series,
(\ref{eq:Vi1}) shows that the spaces $\Ra{\bfV_{i,1}}$ will
be close to each other for a large range of parameter
vectors $\bfsfp_i$, since they are all close to
$\Ra{( \bfA_0 + h^2 \m_1 \bfI)^{-1}\bfB }$).
Regarding the rate of convergence of the series, notice that
the large eigenvalues of
$( \bfA_0 + h^2 \m_1 \bfI)^{-1}$
correspond to globally smooth, low frequency
eigenvectors, whereas the nonzero columns of $\bfI_2$ are
Cartesian basis vectors.
Hence, the series converges rapidly.
Considering more generally the matrices
$\bfV_{i,j} = (\frac{\imath \o_j}{ \nu} \bfE  + \bfA(\bfsfp_i))^{-1}\bfB$,
where $\o_j / \nu$ is small for $j = 2, 3, \ldots, n_{\o}$,
a similar argument holds.
An analogous discussion holds for the $\bfW$ basis.

We demonstrate this experimentally for two model problems:  Cup Image (shown in Figure \ref{fig:cup-recon}(a))
and Amoeba Image (shown in Figure \ref{fig:amoebe-recon}(b)).
Let $\bfV_1, \bfV_2, \ldots$ be the right projection spaces for
the parameter vectors $\bfsfp_1, \bfsfp_2, \ldots$ obtained in steps
$1$ (initial guess), $2$, and so on, of the nonlinear least squares
optimization, that is,
\[
  \bfV_k =
  \left[ \left(\imath \frac{\omega_1}{\nu}\bfE + \bfA(\bfsfp_k)\right)^{-1}\bfB
  \, , \, \ldots \, , \, \left(\imath \frac{\omega_{n_\omega}}{\nu}\bfE +
  \bfA(\bfsfp_k)\right)^{-1}\bfB \right].
\]
In Figure~\ref{fig:MinSinVal}, we show how close these (right) projection
spaces remain to the initial space $\bfV_1$ for two
separate test problems.
We give the sine of the largest canonical angle (the subspace gap)
at each optimization step. As can be seen, for both problems the sine
of the largest canonical angle remains quite small.
Figure~\ref{fig:OptStepImgsCup} shows the changes of the absorption image over the
first six distinct approximate solutions, $\bfsfp_k$.
The first image corresponds to the initial guess for
the parameter vector, $\bfsfp_1$.
Note how drastic changes in the image correspond
to quite small changes in the (right) projection spaces.
\begin{figure}
\bmp{[ht]}{2.5in}
\includegraphics[width=5in]{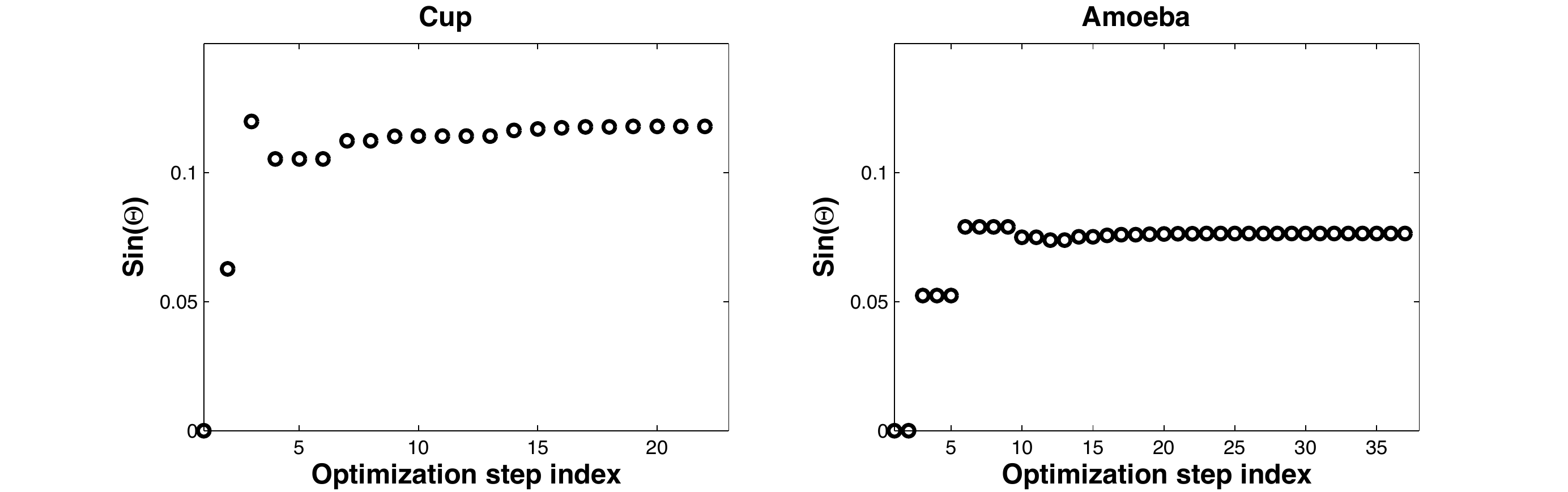}
\emp
\caption{Evolution of the subspace gap (sine of the largest canonical
angle $\Theta$) between initial and subsequent reduction spaces
over the course of the optimization. \label{fig:MinSinVal}}
\end{figure}

\begin{figure}[th]
\centering
\includegraphics[width=4.25in]{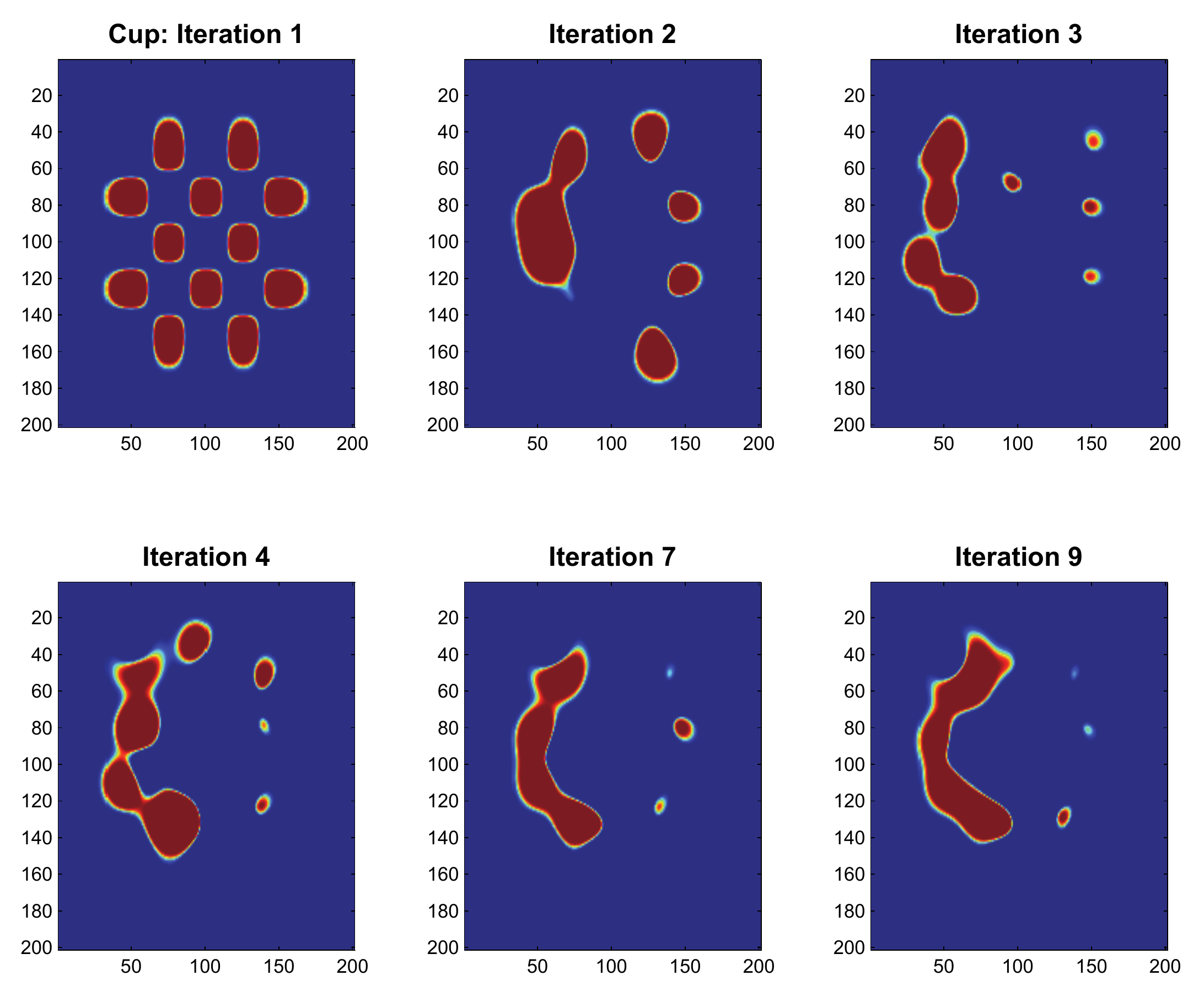}
\caption{Initial image (1) and the first five distinct reconstructions.
\label{fig:OptStepImgsCup}}
\end{figure}

Next, we provide an estimate of the accuracy of the transfer function
based directly on such canonical angles.
For the discussion below, we assume for simplicity that the
sines of the maximum canonical angles for the $\bfW_k$ spaces are similar
to those for the $\bfV_k$.
Now consider using, at step $k$, the reduced model computed from
$\bfV_j$ and $\bfW_j$ obtained from $\bfsfp_j$ in \eqref{eq:PGProj}
for the frequency interpolation points $\omega_1, \ldots, \omega_{n_\o}$,
$\wh{\bfPsi}_j(\omega;\bfsfp)$, and the resulting error at $(\omega_i,\bfsfp_k)$,
$\bfPsi(\omega_i;\bfsfp_k) - \wh{\bfPsi}_j(\omega_i;\bfsfp_k)$.
Denote the reduced model that would be computed from
$\bfV_k$ and $\bfW_k$ obtained from $\bfsfp_k$ in \eqref{eq:PGProj},
by $\wh{\bfPsi}_k(\omega;\bfsfp)$. By Theorem~\ref{thm:parammor},
$\wh{\bfPsi}_k(\omega_i;\bfsfp_k) = \bfPsi(\omega_i;\bfsfp_k)$; the reduced
model $\wh{\bfPsi}_k(\omega;\bfsfp)$ is exact at $\bfsfp_k$ (for all $\omega_i$).
Furthermore, Theorem~\ref{thm:parammor} states that
$\wh{\bfPsi}_j(\omega_i;\bfsfp_k)$ would be exact
if $\Ra{\bfV_j} = \Ra{\bfV_k}$ and $\Ra{\bfW_j} = \Ra{\bfW_k}$.
This will generally not be the case; however, the computed results in
Figure~\ref{fig:MinSinVal} suggest that these spaces may be
very close.

Theorem~3.3 in \cite{beattie2010isi}
provides a bound on the relative error
in approximating $\wh{\bfPsi}_k(\omega_i;\bfsfp_k)$ by
$\wh{\bfPsi}_j(\omega_i;\bfsfp_k)$
based on the canonical angles between
$\Ra{\bfV_k}$ and $\Ra{\bfV_j}$ and the canonical angles between
$\Ra{\bfW_k}$ and $\Ra{\bfW_j}$.
In addition, using that $\wh{\bfPsi}_k(\omega_i;\bfsfp_k) = \bfPsi(\omega_i;\bfsfp_k)$
for the frequency interpolation points $\omega_1, \ldots, \omega_{n_\o}$,
we obtain the following result.
For $i=1,\ldots,n_\omega$,
\eqs \nonumber
  \frac{ \| {\bfPsi}(\omega_i;\bfsfp_k) - \wh{\bfPsi}_j(\omega_i;\bfsfp_k) \|_2 }
       {\frac{1}{2}\left( \| \wh{\bfPsi}_k(\,\cdot\,;\bfsfp_k) \|_{\Hinf} +
              \| \wh{\bfPsi}_j(\,\cdot\,;\bfsfp_k) \|_{\Hinf} \right)} & \leq &
  M \max\left( \sin\V(\bfV_k,\bfV_j), \sin\V(\bfW_k,\bfW_j) \right) ,
\\ \label{eq:AngleForwBnd}
  &&
\eqe
where $M$ depends on the conditioning of the reduced order model
with respect to $\bfB$ and $\bfC$ and the angles between some of
the spaces associated with the reduced order model (see \cite{beattie2010isi}),
and $\| \wh{\bfPsi}_i(\,\cdot\,;\bfsfp_k) \|_\Hinf = \sup_{\omega \in \IR} \| \wh{\bfPsi}_i(\omega;\bfsfp_k)\|_2$
is the $\Hinf$ norm.

Using \eqref{eq:AngleForwBnd}, the results given in Figure~\ref{fig:MinSinVal},
and (assuming) similar bounds for the spaces $\bfW_k$, we can
assess the errors from using a reduced model computed from $\bfV_1$ for
each step in the optimization
(with parameter vectors $\bfsfp_1, \bfsfp_2, \ldots$). The
usefulness, of course, depends on the condition number $M$.
To test how descriptive the subspace gaps are,
we also compute the relative interpolation errors given in the left-hand
side of \eqref{eq:AngleForwBnd} for the Cup and Amoeba image reconstructions,
and we compare these with the subspace gaps given  Figure \ref{fig:MinSinVal}.
For both test cases, the interpolation errors, shown in Figure \ref{fig:InterpolationError},
follow a pattern very similar to that of the subspace gaps.
These numbers suggest the possibility of a modest condition number $M$;
thus illustrating that, at least for these examples, the gaps could be
a good indicator for the behavior of the interpolation error.

Note that the interpolation errors in our numerical experiments with
reduced models will be much smaller than what is shown here.
Here, we report results as if only the $\bfV$ subspace is used in the reduction step, and the adjoint information due to $\bfW$ is ignored. Moreover, for our numerical experiments,
the subspaces $\bfV$ and $\bfW$ will be constructed by sampling more than one parameter value
(up to five). Therefore, they carry much more global information than
is used in our simple analysis of the subspace gaps, here.
The effectiveness of interpolatory reduced-models for the inversion problem is illustrated in Section~\ref{sec:NumExp} using four numerical examples.

\begin{figure}
\centering
\begin{tabular}{cc}
\includegraphics[width=2.0in]{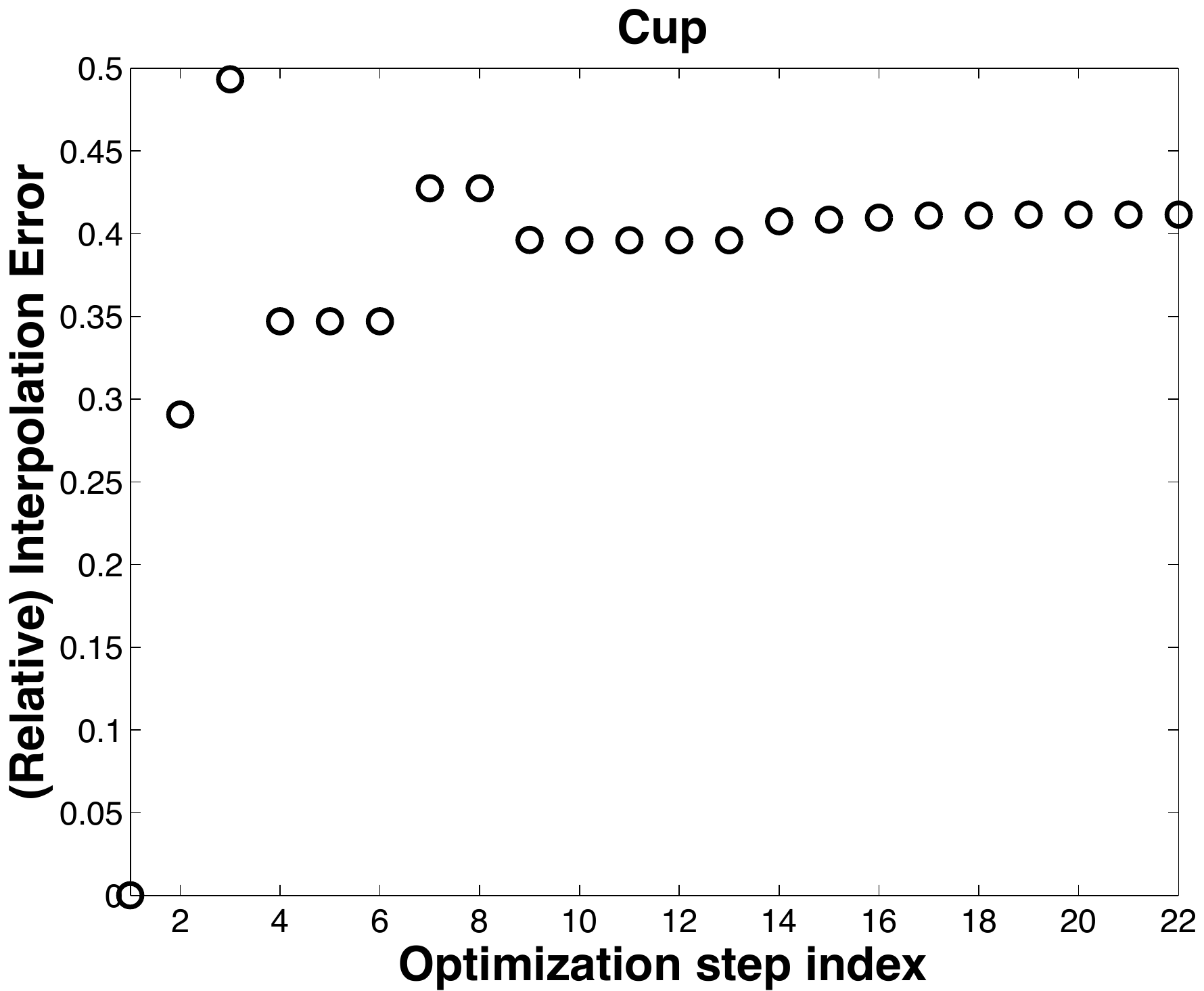} &
\includegraphics[width=2.0in]{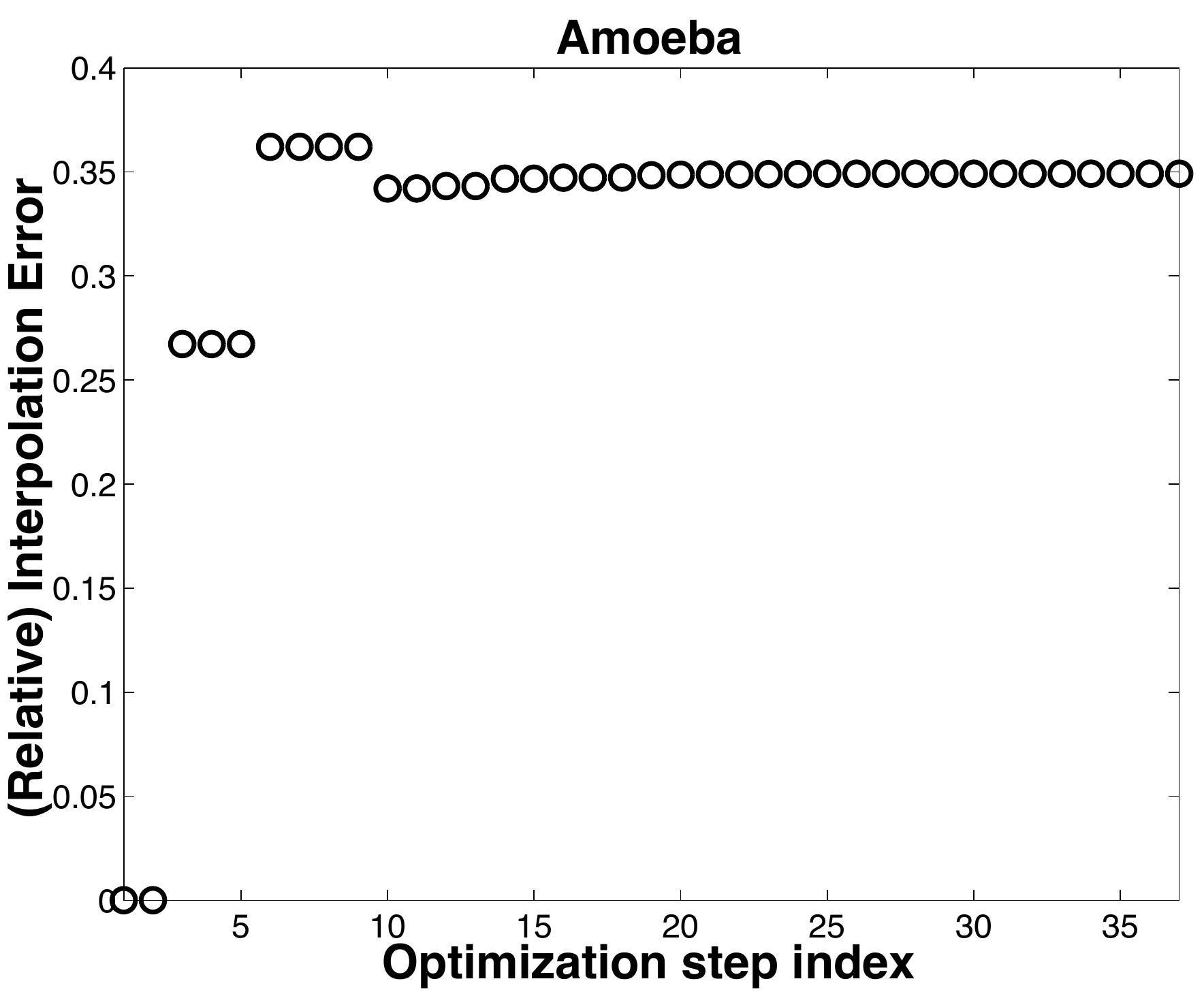}
\end{tabular}
\caption{Evolution of the relative interpolation error
over the course of the optimization. 
Note that larger model reduction spaces carrying more global information are used in the actual numerical
experiments in the next section, leading to much smaller
interpolation errors.
\label{fig:InterpolationError}}
\end{figure}

The results of our analysis 
suggest that using  parameter vectors from the first few
steps of the optimization to compute global bases for a reduced order model
will lead to a reduced model that can be used effectively throughout the optimization and
likely also for distinct image reconstructions as well (as the initial parameter vectors
would all be close to the initial guess). The numerical experiments presented in the
next section 
suggest this as well.

%

\subsection{Implementation and Cost Overview}\label{sec:Impl-n-Cost}
As in most parametric model reduction approaches, the offline phase of our
method consists of sampling the parameter space and constructing the model
reduction bases $\bfV$ and $\bfW$ corresponding to the sampled parameter
points. Let $\{\bfpi_i\}_{i=1}^K$ denote the parameter samples. In our
numerical examples, as explained in Section \ref{sec:NumExp} below, we use
the leading $K \leq 5$ iterations
of the full optimization algorithm to obtain the sample points.
However, the samples can come from any effective sampling technique.
Once the  parameter samples are determined, the local bases $\bfV_i$, for
$i=1,\ldots,K$, are constructed  as in (\ref{eq:localVi}),
and a basis $\widetilde \bfV$
is obtained by concatenating these local bases, i.e.
$\widetilde\bfV = [\bfV_1, \bfV_2,\ldots, \bfV_{K}]$.
However, in most cases, $\widetilde \bfV$ will contain linearly dependent columns.
Thus, the final global basis $\bfV$ is taken as the leading left singular
vectors of $\widetilde \bfV$ corresponding to the singular values bigger than
a given tolerance. If the SVD removes small, but  non-zero,
singular values below the given tolerance, then the interpolation at the sampled frequency and parameter points  will be approximate.
The same steps are repeated to construct $\bfW$. Construction
of $\bfV$ and $\bfW$ concludes the offline stage.

In the online stage, we use the parametric reduced model to replace
the expensive large-scale function and Jacobian evaluations.  Using
$\bfV$ and $\bfW$, for the current iterate $\bfsfp_k$ of the optimization
problem, we need $\widehat{\bfE} = \bfW^{T}\bfE\bfV$,
$\widehat{\bfB}= \bfW^{T}\bfB$, $\widehat{\bfC} = \bfC\bfV$, which
are constant matrices and are computed only once, and
$\widehat{\bfA}(\bfsfp_k) = \bfW^{T} \bfA(\bfsfp_k) \bfV$, which must be updated
for each new iterate $\bfsfp_k$. Then, in
the optimization algorithm, the large-scale function evaluation
$\bfPsi(\bfsfp_k)$ and (if necessary) the Jacobian evaluation
$\nabla_{\bfsfp}{\bfPsi}(\bfsfp_k)$ are replaced with reduced-model
evaluations $\widehat{\bfPsi}(\bfsfp_k)$ and
$\nabla_{\bfsfp}\widehat{\bfPsi}(\bfsfp_k)$, respectively.

The major computational cost results from constructing the model
reduction bases $\bfV$ and $\bfW$. Assuming $n_\omega$ frequency
interpolation points, $n_{src}$ sources, $n_{det}$ detectors, and
$K$ parameter samples, constructing $\bfV$ and $\bfW$ requires solving
$K n_\omega n_{src} + K n_\omega n_{det}$
large, sparse, $n \times n$, linear systems.
In our examples, $n_{det} = n_{src}$, so the total
number of large linear systems to be solved equals
$2K n_\omega n_{src}$.
As explained above, the cost of the online stage is modest,
since it only requires  $r \times r $ linear solves and the multiplication
$\bfW^{T} \bfA(\bfsfp) \bfV$ can be performed at a cost independent of $n$;
thus the main computational complexity is due to solving
$2K n_\omega n_{src}$
large linear systems. On the other hand, the optimization using
the full model
requires solving  $K_{fun} n_\o n_{src} + K_{Jac} n_\o n_{det}$ large linear systems,
where $K_{fun}$ is the number of objective function evaluations and
$K_{Jac}$ is the number of computations of the Jacobian over the entire optimization.
In our numerical examples, the ratio $(K_{fun} + K_{Jac}) / 2K$ is always larger
than three; thus the full parameter inversion requires at least three times
more large-scale linear solves than the reduced-model based parameter
inversion does; thus the offline costs are quickly amortized. Indeed, as
we explain in Section \ref{sec:NumExp} below, the computational gains
are even much higher, since we are able to recycle the model
reduction bases for reconstructing different images; in other words,
there are no offline costs for subsequent reduced-model based
parameter inversions for the same mesh
with the same source and detector locations.

\section{Numerical Experiments}\label{sec:NumExp}

The purpose of our numerical experiments is to provide a proof of concept of our approach
and some analysis and insight. We discuss four distinct reconstructions in two groups of two.
Within each group, the same mesh, the same number of basis functions, and the same sources and detector locations are used. However, the two original images to be (approximately) reconstructed within each group are very different. Apart from showing that the reconstruction using the reduced model is very close to that using the full model, we demonstrate that the second reconstruction can be done using the same projection bases as computed for the first reconstruction, supporting  the discussion of section \ref{sec:anal_gb}.
{\em This means that no further solves with large systems are required, i.e., without any
further off-line cost.} Of course, some computational effort is required computing
the reduced model, however these are a modest number of additional
matrix-vector products with a diagonal matrix and a modest number of dot-products.
See sections~\ref{ssec:GlobalBases} and \ref{sec:Impl-n-Cost} for a detailed
discussion of the cost. These examples strongly suggest that after good bases
have been computed, they can be used for many distinct reconstructions.
This would be of significant practical importance, as it may make the cost of computing
the bases, which does require a number of iterative solves with large
matrices, more or less insignificant.

The experiments are set up as follows. We generate synthetic data by
first constructing one or more anomalies in a standard $0-1$ pixel basis.
Then, this $0-1$ image is mapped to an image having pixel-wise
the absorption value $\m_{out}$ for healthy tissue (for the $0$ pixels)
and the absorption value $\m_{in}$ for anomalous tissue (for the $1$ pixels).
Next, the absorption values inside the anomaly are given a small normally
distributed random variation. This pixel-based image of absorption is
then used for solving the forward problem for each source and frequency
and computing the measured values. These measurements are further perturbed
by adding white noise (0.1\%). Finally, we reconstruct the shape
of the absorption images using the compactly supported radial basis
functions described in \S\ref{ssec:pals}.
The optimization problem is solved using the TREGS algorithm
\cite{StuKil11c} with both the exact objective function (full model)
and the approximate objective function represented by the reduced order model.

\subsection*{Small Mesh}

The first two problems are solved on a uniform $50 \times 50$ grid,
physically corresponding to a domain of $5\mathrm{cm}\times 5\mathrm{cm}$.
We use standard centered finite differences to discretize the DOT problem (\ref{eq:ModelPDE1}--\ref{eq:ModelPDE4}). We use $15$ CSRBFs to reconstruct the absorption image, leading to $60$ parameters (four per 2D basis function). The model has $24$ sources and $24$ detectors, i.e., $\bfB \in \IR^{2500 \times 24}$ and $\bfC \in \IR^{24 \times 2500}$.

To compute a reduced model for the first problem, we choose
the parameter vectors, $\bfsfp_1$ and $\bfsfp_2$, from the first two iterations
of the optimization using the exact objective function.\footnote{This is a
reasonable approach in practice, as it allows us to switch to
the reduced model in the optimization after a few steps of optimization
using the full model. Alternatively, one can choose a number of parameter
interpolation points a priori. A brief discussion of approaches is given
in section \ref{sec:Impl-n-Cost}.}
Next we solve, for each parameter point and one frequency,
for the $24$ right hand sides given by the sources,
$\bfV_{i,1} = (\frac{\imath \,\!\omega_1}{\n}\bfE + \bfA(\bfsfp_i))^{-1}\bfB$,
and the $24$ right hand
sides given by the detectors,
$\bfW_{i,1} = (\frac{\imath \,\!\omega_1}{\n}\bfE + \bfA(\bfsfp_i))^{-T}\bfC^T$,
for $i=1,2$. To maintain the symmetry of the full model in the reduced model,
we employ one-sided projection, i.e., $\bfW = \bfV$. We achieve this by
first combining the local basis matrices together,
$\widetilde{\bfV} = [\bfV_1 \, \bfV_2  \, \bfW_1 \, \bfW_2] \in \IR^{2500 \times 96}$,
and then reducing the dimension using a rank-revealing decomposition of
$\widetilde{\bfV}$ with some modest tolerance (here, we use an SVD, but
cheaper approaches will be used for larger systems). This results
in a final reduction basis $\bfW = \bfV \in \IR^{2500 \times 80}$
and therefore a reduced model of dimension $r=80$. The results of
the reconstruction are given in Figure~\ref{fig:tri-block}.

\begin{figure}
\bmp{[t]}{1.5in}
\includegraphics[width=1.25in,height=1.25in]{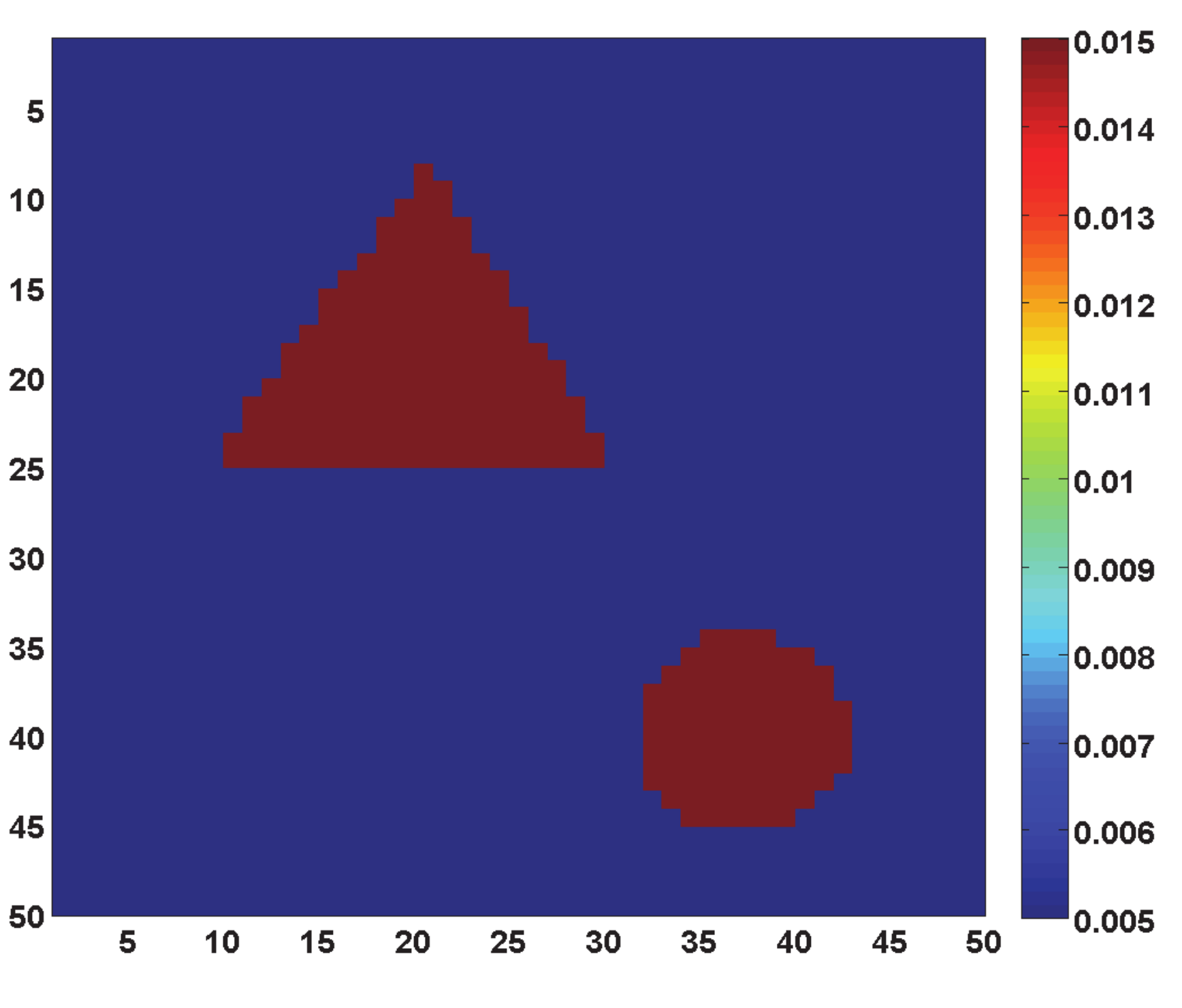}
\\
(a) Original shapes
\emp
\hspace{\stretch{1}}
\bmp{[t]}{1.5in}
\includegraphics[width=1.25in,height=1.25in]{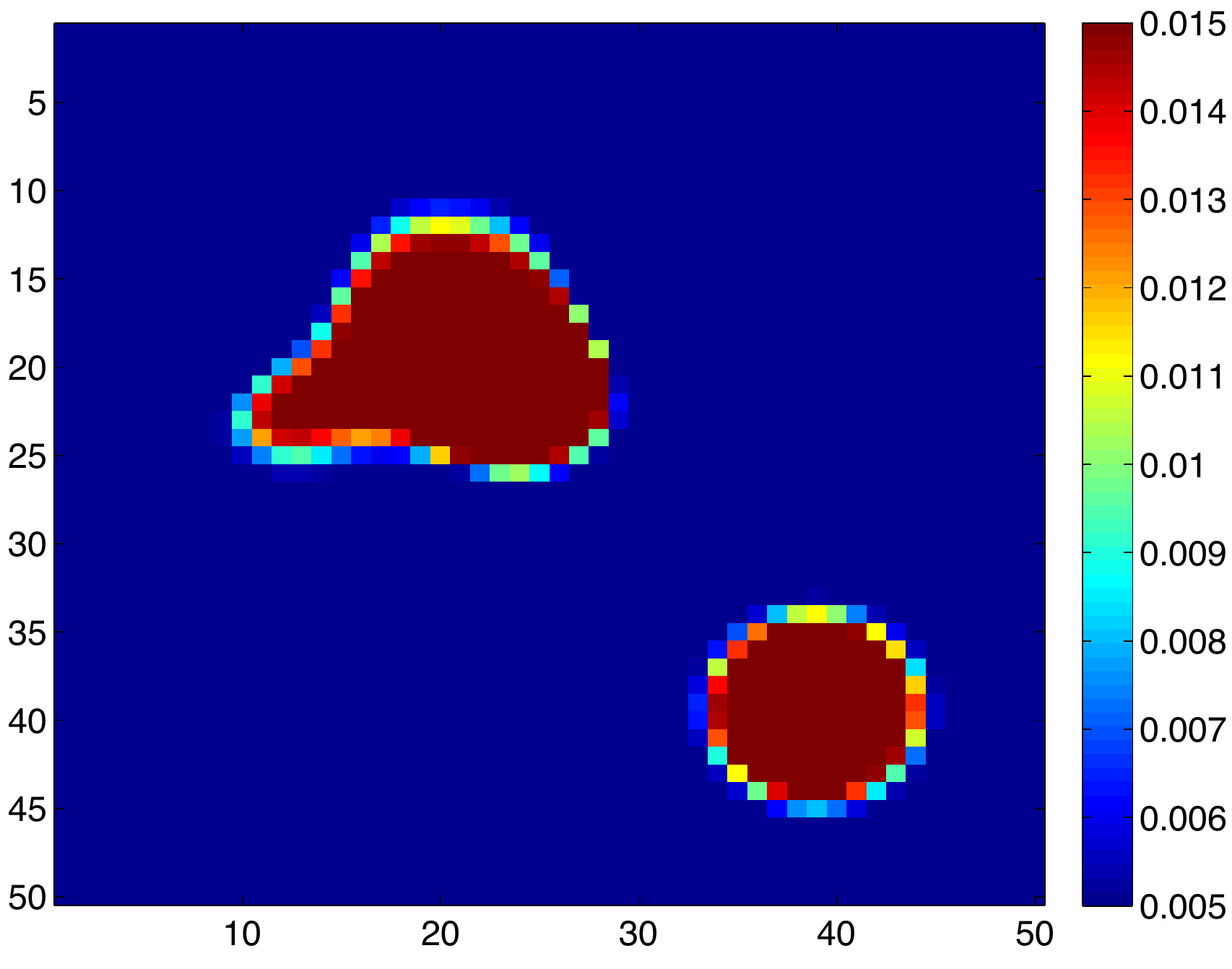}
\\
(b) Reconstruction using the exact objective function (full model)
\emp
\hspace{\stretch{1}}
\bmp{[t]}{1.5in}
\includegraphics[width=1.25in,height=1.25in]{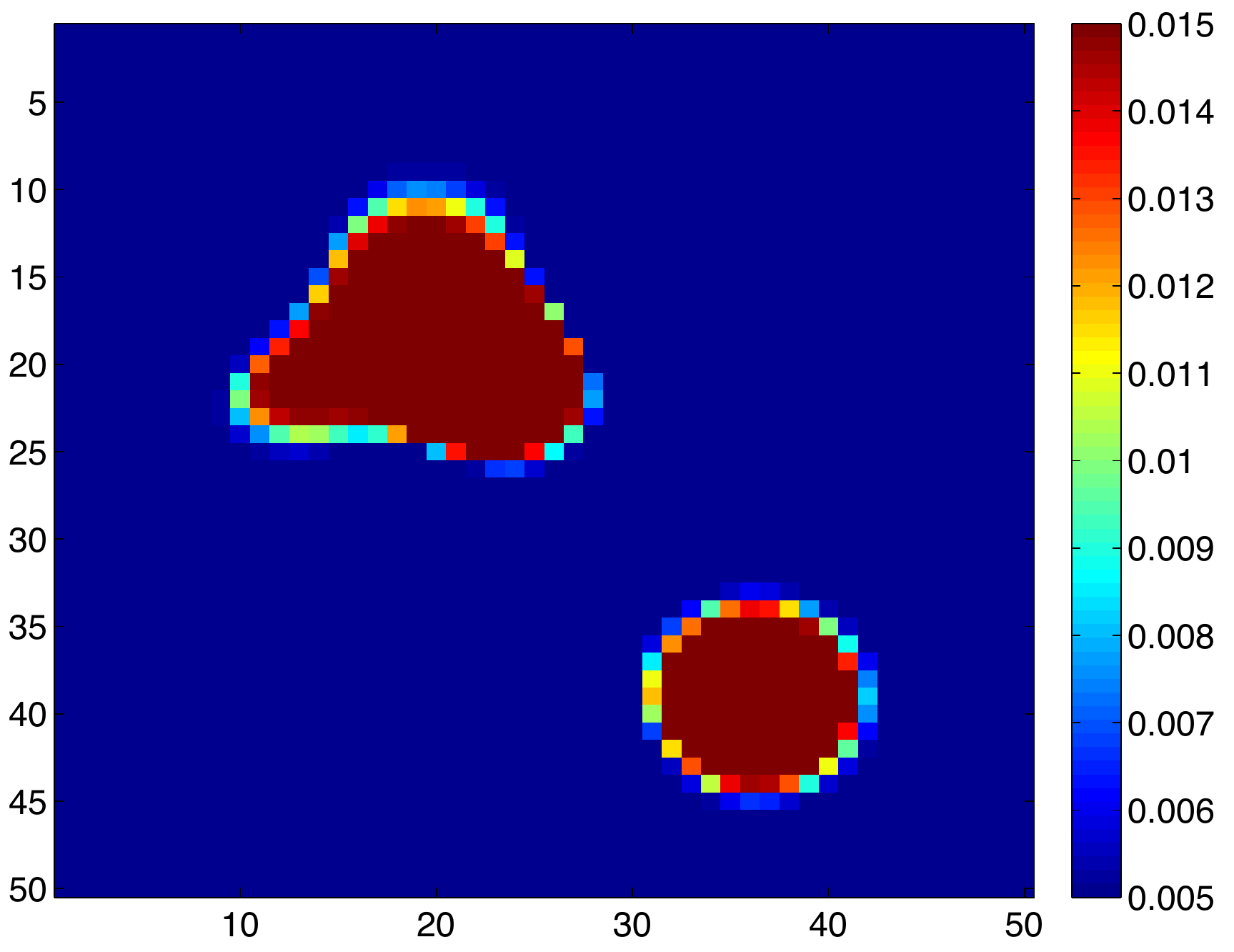}
\\
(c) Reconstruction using an approximate objective function (reduced order model)
\emp

\caption{Reconstruction of a simple test anomaly on a small mesh
with $24$ sources and $24$ detectors}
\label{fig:tri-block}
\end{figure}

The reconstructions using the exact objective function and the approximate objective function (reduced
model) converge to the same tolerance, and are very similar.
Both reconstructions approximate the original shape well.
Note that an exact or very accurate reconstruction is not possible,
as we are not reconstructing in the pixel basis, and our model does not
reconstruct the small heterogeneity in the anomaly. The result clearly suggests that
reconstructions using a reduced model can obtain a quality similar to those using the full model.
This means that the bases $\bfV$ and $\bfW$, constructed using the
two parameter vectors from the first two iterations, have provided accurate interpolations
and thus accurate function and Jacobian evaluations for the rest of the optimization
validating the discussion of section \ref{sec:anal_gb}.
Using the full model, the optimization algorithm solves $768$ linear
systems of dimension $2500$ to reconstruct the absorption
image. Using the reduced model, the optimization algorithm solves only $96$ linear systems of dimension
$2500$ to compute the projection basis for the reduced model and $816$ linear systems of dimension $80$ to reconstruct
the absorption image. Hence, the number of `full size' linear systems that must be solved is drastically reduced. For
this problem the full system size is small anyway, but for large matrices (a large forward problem) the computational
savings will be significant. This will be all the more so, if many sources and detectors are used and multiple
frequencies.

Next, we investigate a reconstruction problem at the same discretization level (i.e. $n=2500$)
and with the same source and detector locations as before,
but for a very different image, as shown in Figure \ref{fig:triple-clubs}-(a).
Since our experiments and analysis in \S~\ref{sec:anal_gb} suggest that
the projection bases vary only modestly with $\bfsfp$,
for this second problem, we do not compute new model reduction bases; instead,
we use the same basis $\bfV = \bfW$ from the previous example,
even though the image to reconstruct is quite different (in the pixel basis).
Thus, there is {\it no offline cost} for constructing the model reduction basis for the second
image reconstruction. However, we want to clarify that recycling the model reduction basis
does not mean that we are using the same parametric reduced model.
As the optimization goes through different parameter iterates, $\bfsfp_k$,
the matrices $\bfA(\bfsfp_k)$ will differ from those arising for the first image;
therefore, the reduced quantity $\bfA_r(\bfsfp) = \bfW^T \bfA(\bfsfp) \bfV$, and thus
the reduced parametric model, will be different.
As can be seen from Figure~\ref{fig:triple-clubs}, the use of a reduced model
as an approximate objective function yields a reconstruction that is as
good as that using the exact objective function, even
using a reduced system basis that was computed for a substantially different reconstruction.
As we use a precomputed basis (from the previous problem), no solves (direct or iterative)
with the full matrix are needed.
Only small (reduced system size) linear systems have to be solved. Of course, for this (small) problem size, the difference is not so large, but for the next two problems of substantially larger size the difference in computational cost is very large.
For {\sc Matlab} code, it is the difference between an iterative solver for a moderately large system versus a quick solve
using `backslash'.

\begin{figure}
\bmp{[t]}{1.5in}
\includegraphics[width=1.25in,height=1.25in]{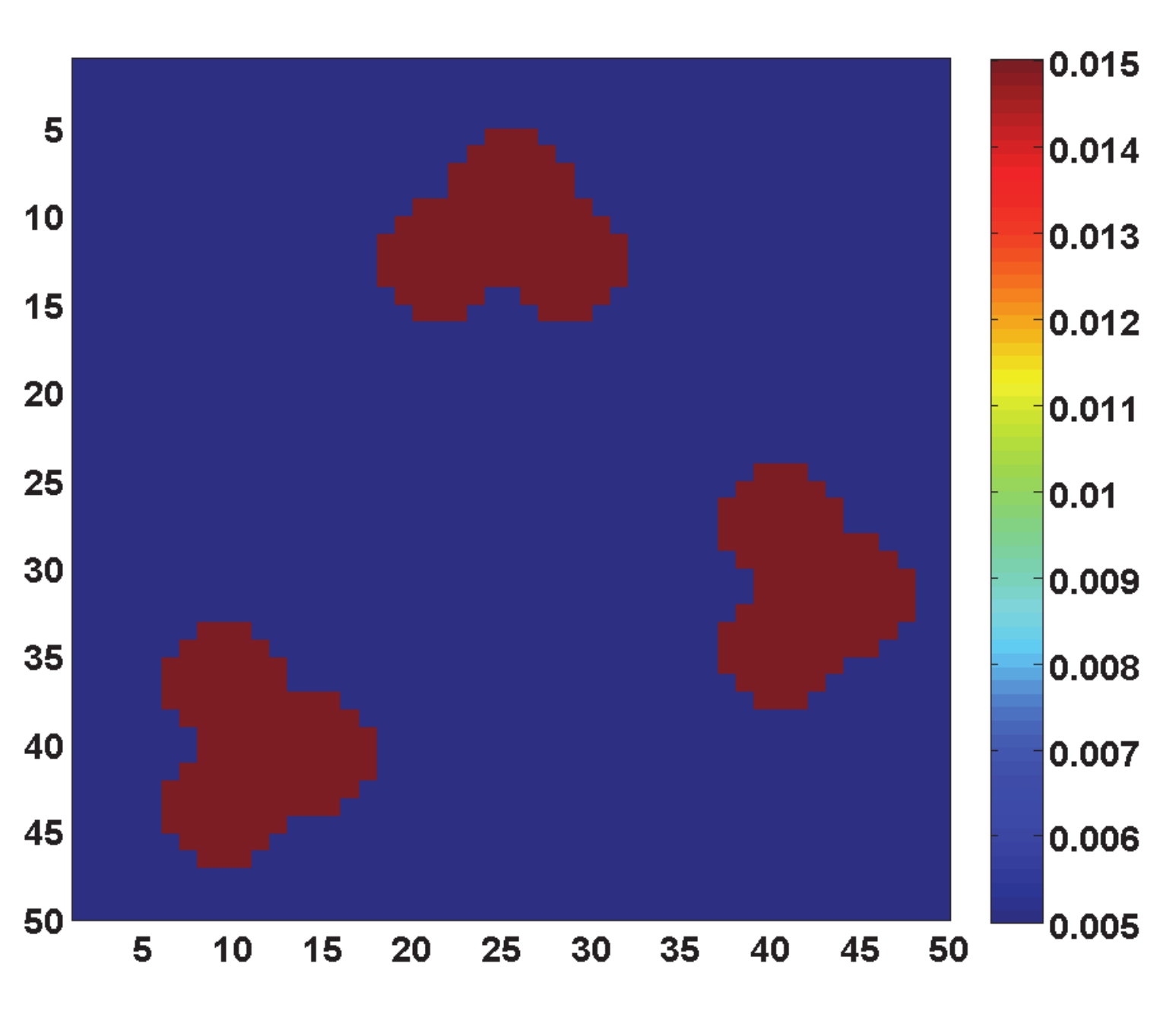}
\\
(a) Original shapes
\emp
\hspace{\stretch{1}}
\bmp{[t]}{1.5in}
\includegraphics[width=1.25in,height=1.25in]{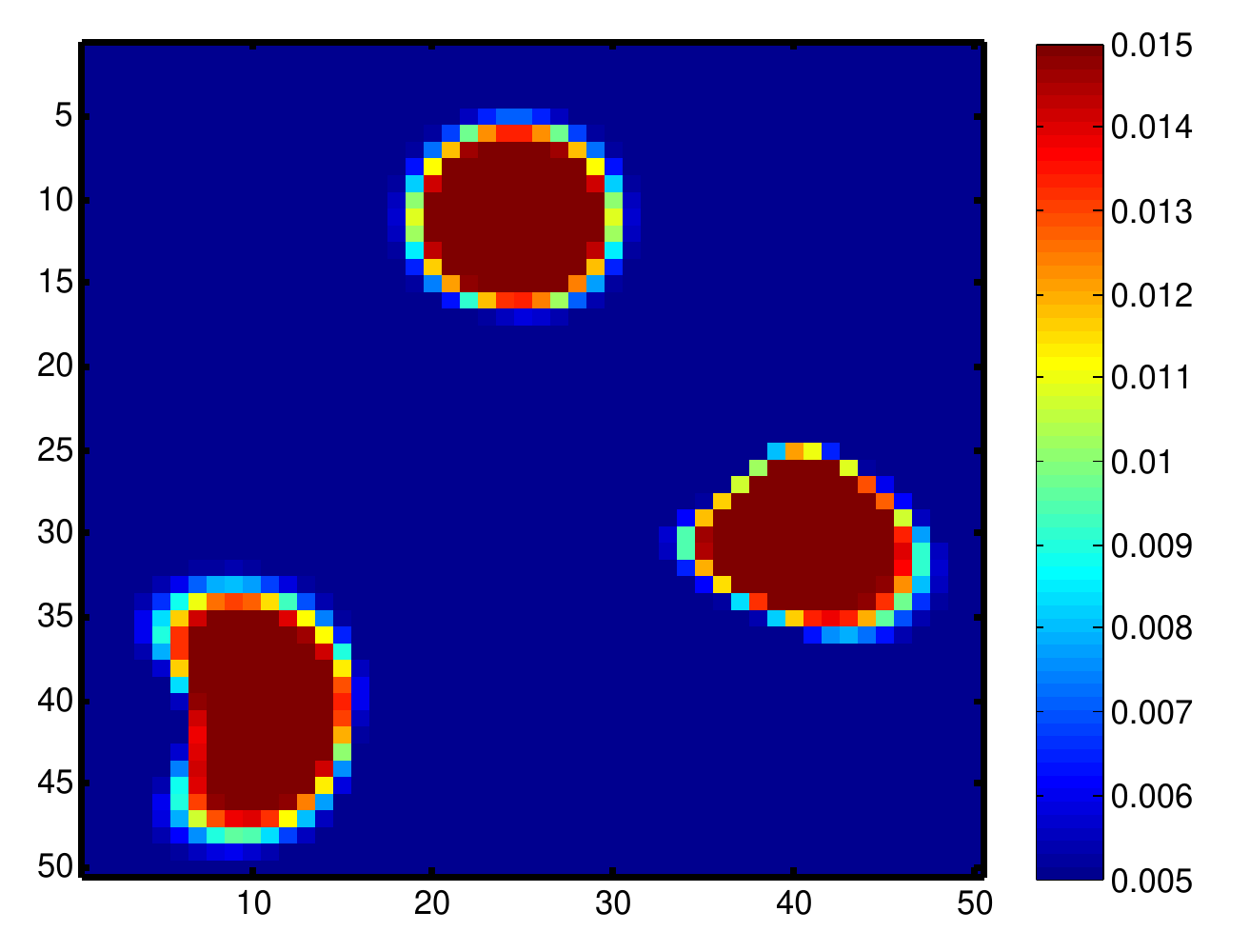}
\\
(b) Reconstruction using the exact objective function (full model)
\emp
\hspace{\stretch{1}}
\bmp{[t]}{1.5in}
\includegraphics[width=1.25in,height=1.25in]{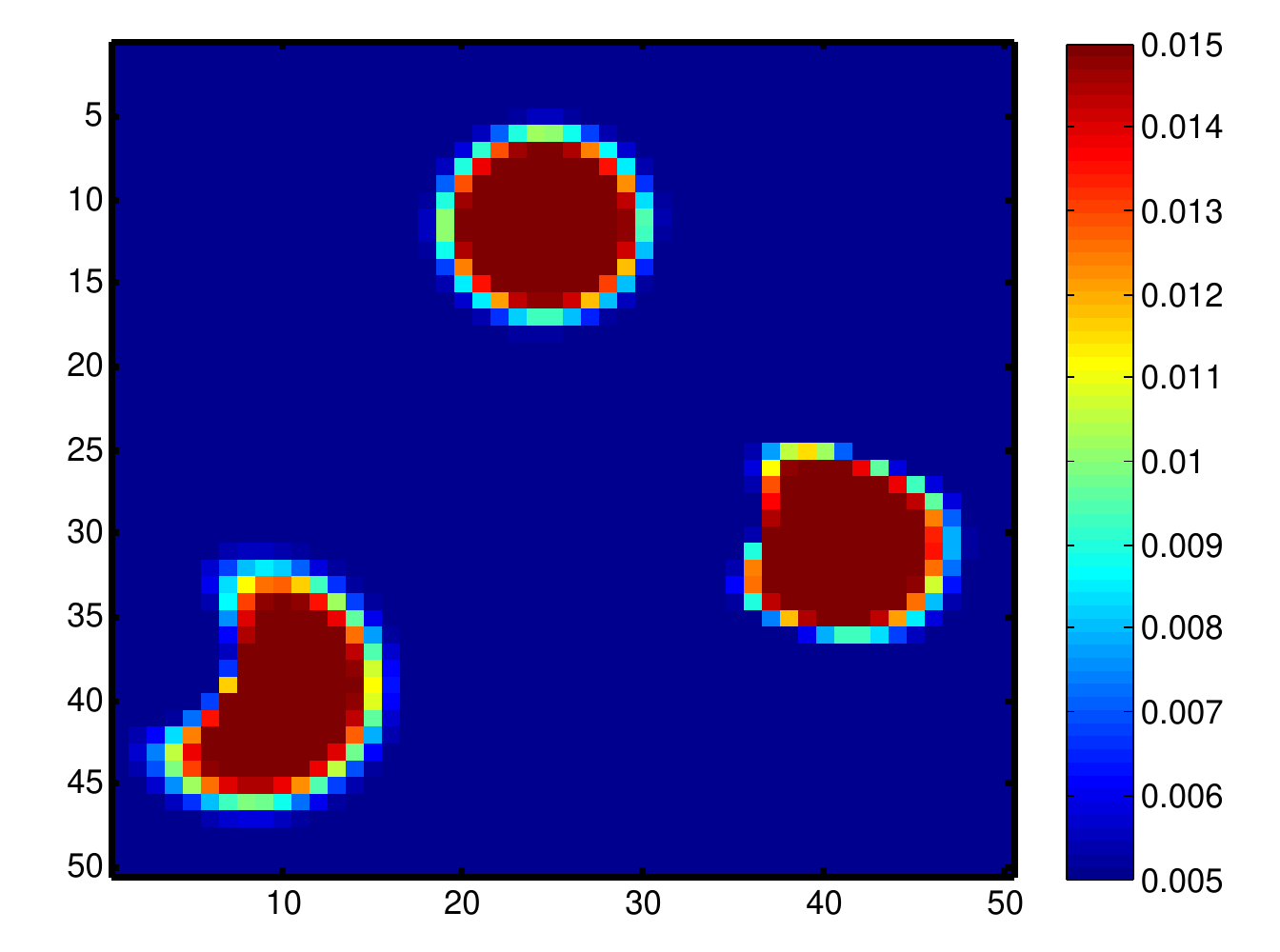}
\\
(c) Reconstruction using an approximate objective function (reduced order model)
\emp

\caption{Results for Problem~2. Reconstruction of a simple test anomaly on a small mesh with $24$ sources and $24$ detectors using the same bases computed for Problem~1. So, {\em only} the reduced model was used.
The discretization and sources and detector locations are the same as for Problem~1.}
\label{fig:triple-clubs}
\end{figure}

Using the exact objective function (full model), the algorithm solved $528$ full size
linear systems (of size $2500 \times 2500$).
Using the reduced model, the algorithm did not solve any full size system and only solved $576$ reduced size systems (of size $80 \times 80$).

\subsection*{Large Mesh}

\begin{figure}
\centering
  \includegraphics[width=1.25in,height=1.25in]{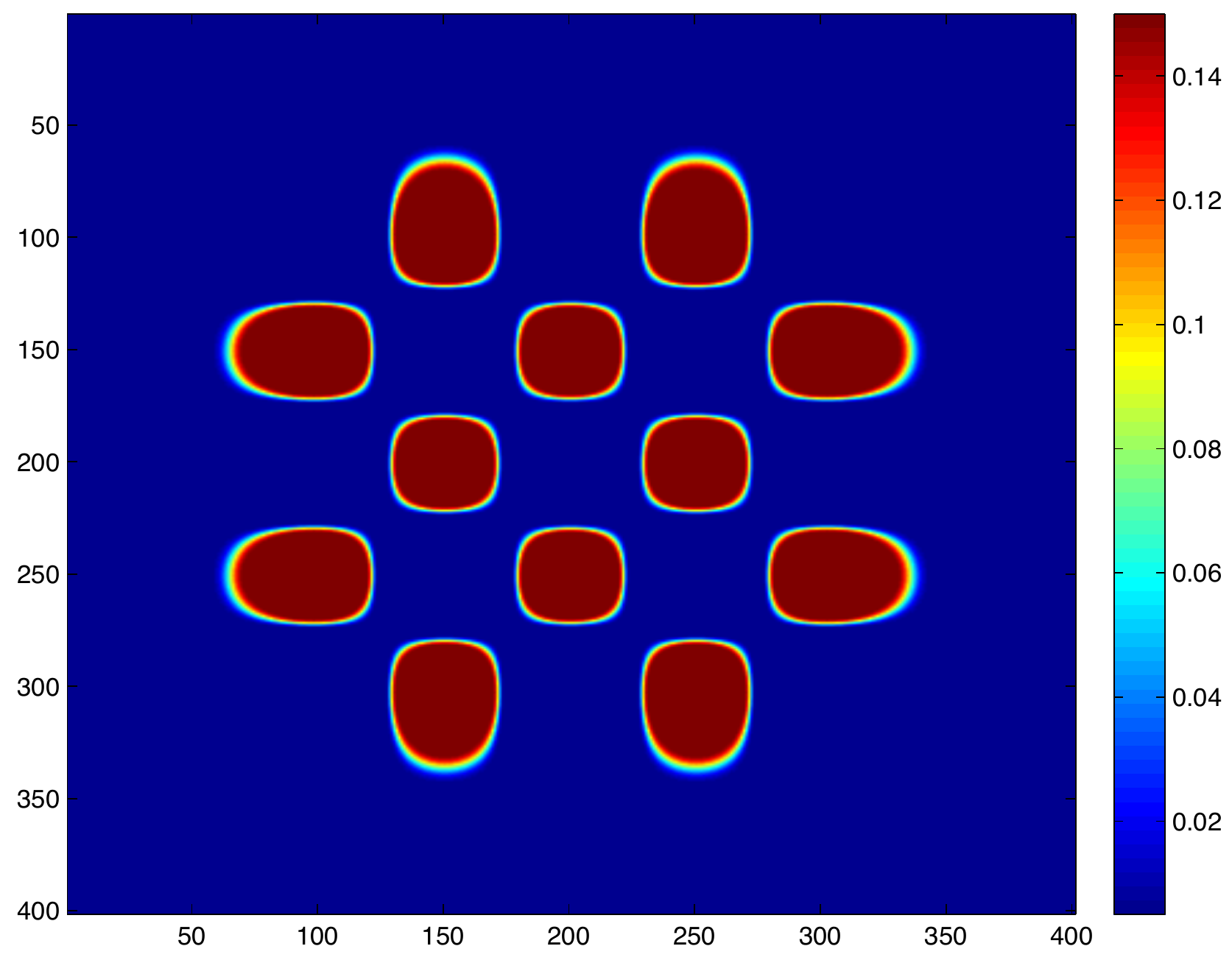}\\
  \caption{Initial configuration with $25$ basis functions arranged in a $5 \times 5$ grid
  with alternatingly negative and positive $\alpha$'s. }
  \label{fig:initial}
\end{figure}

\begin{figure}
\bmp{[t]}{1.5in}
\includegraphics[width=1.25in,height=1.25in]{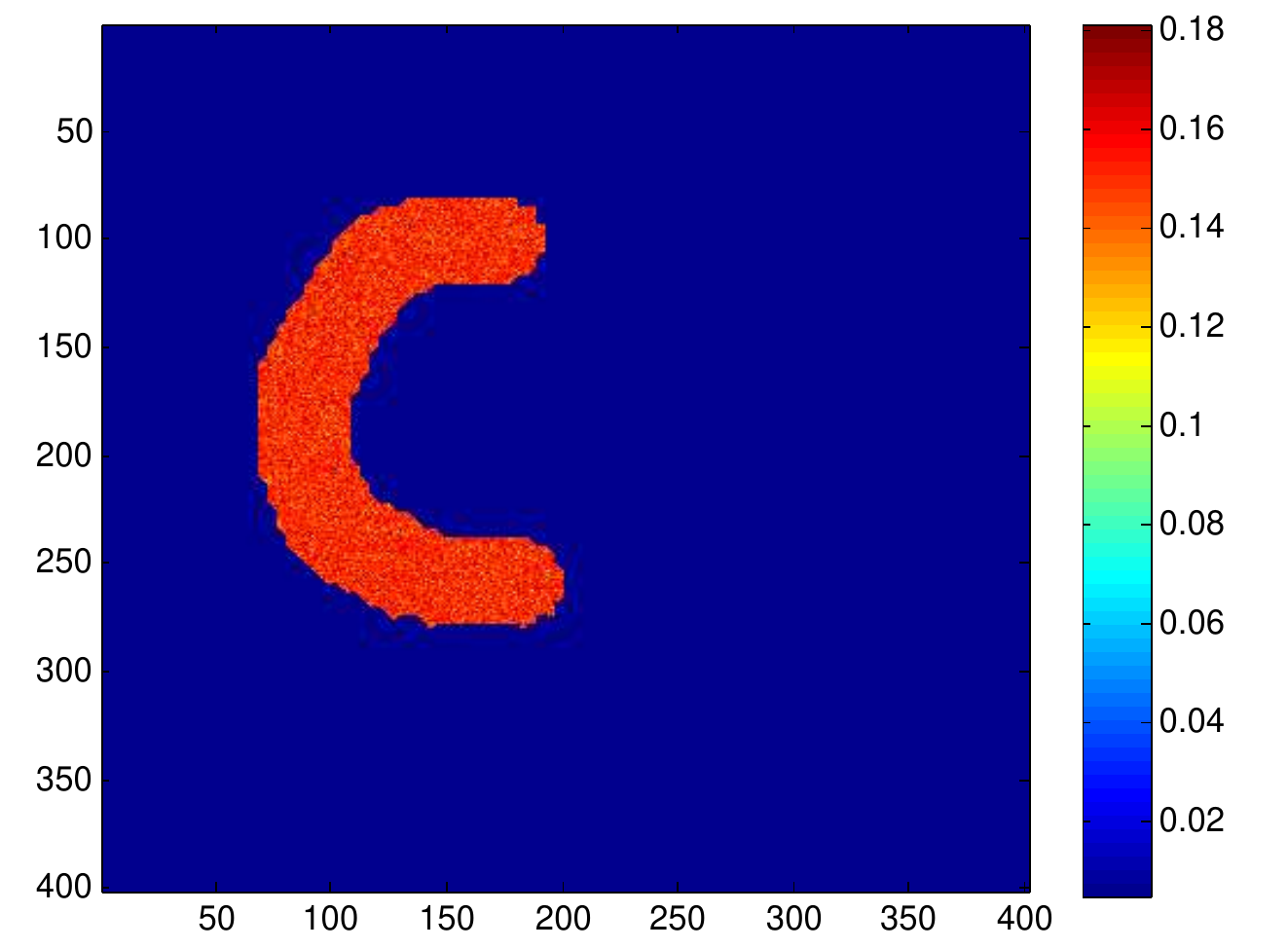}
\\
(a) Original shape.
\emp
\hspace{\stretch{1}}
\bmp{[t]}{1.5in}
\includegraphics[width=1.25in,height=1.25in]{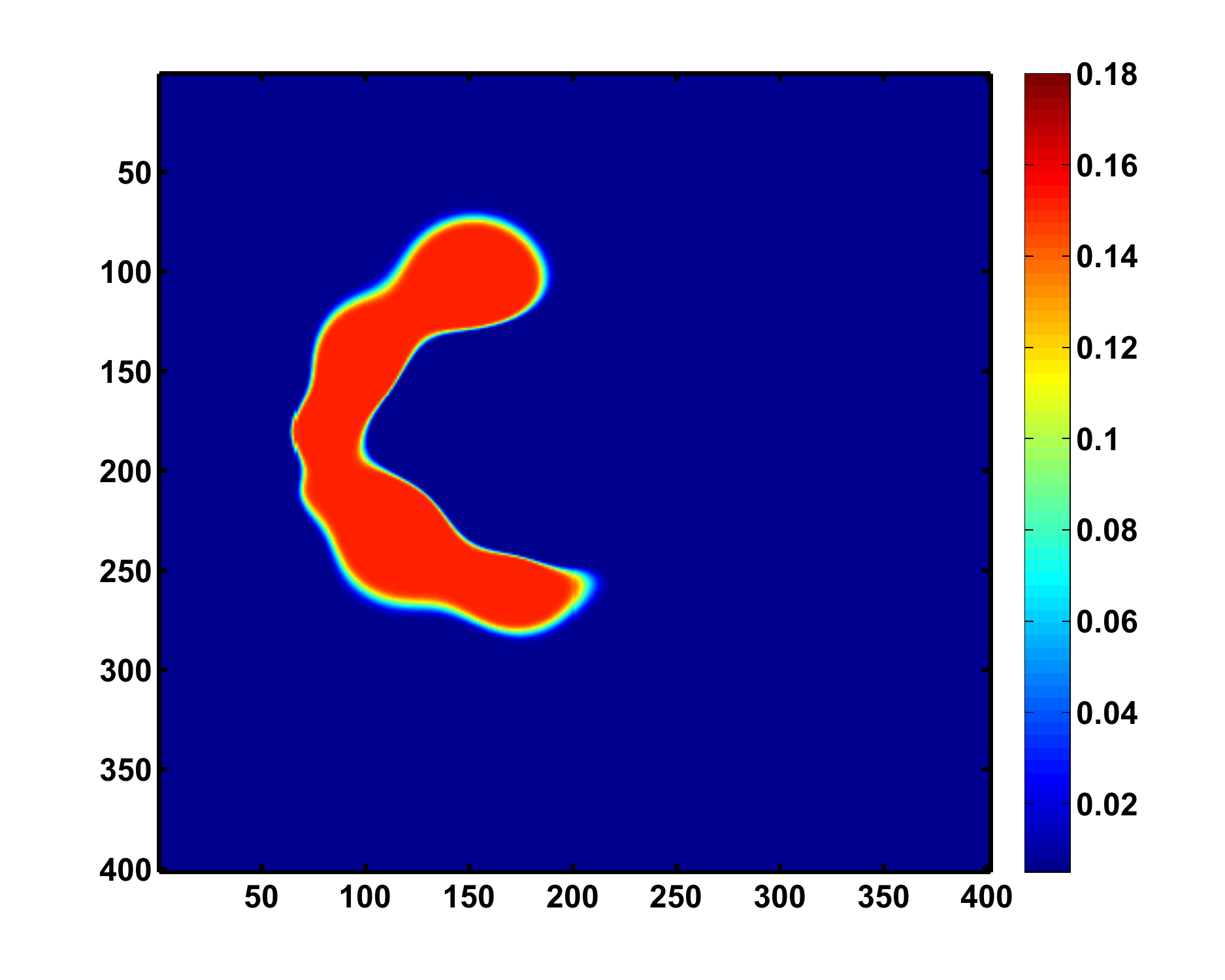}
\\
(b) Reconstruction using the exact objective function (full model).
\emp
\hspace{\stretch{1}}
\bmp{[t]}{1.5in}
\includegraphics[width=1.25in,height=1.25in]{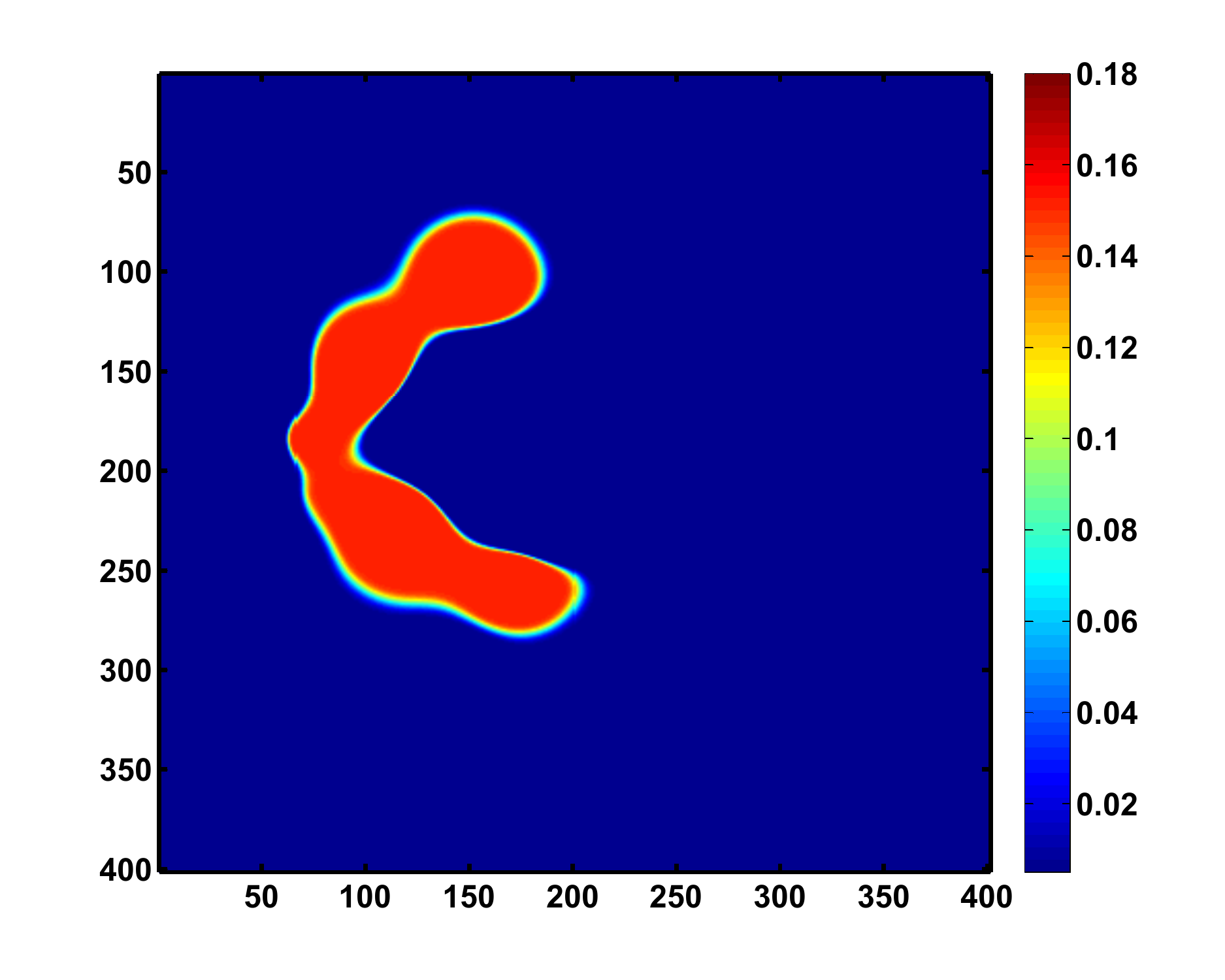}
\\
(c) Reconstruction using an approximate objective function (reduced order model).
\emp

\caption{Results for Problem~3. Reconstruction of a test anomaly on a $401 \times 401$ mesh, resulting in $160801$ degrees of freedom in the forward model, with $36$ sources and detectors, and $25$ basis functions. The reduced model has $250$ degrees of freedom for the forward model.}
\label{fig:cup-recon}

\end{figure}

\begin{figure}
\bmp{[t]}{1.5in}
\includegraphics[width=1.25in,height=1.25in]{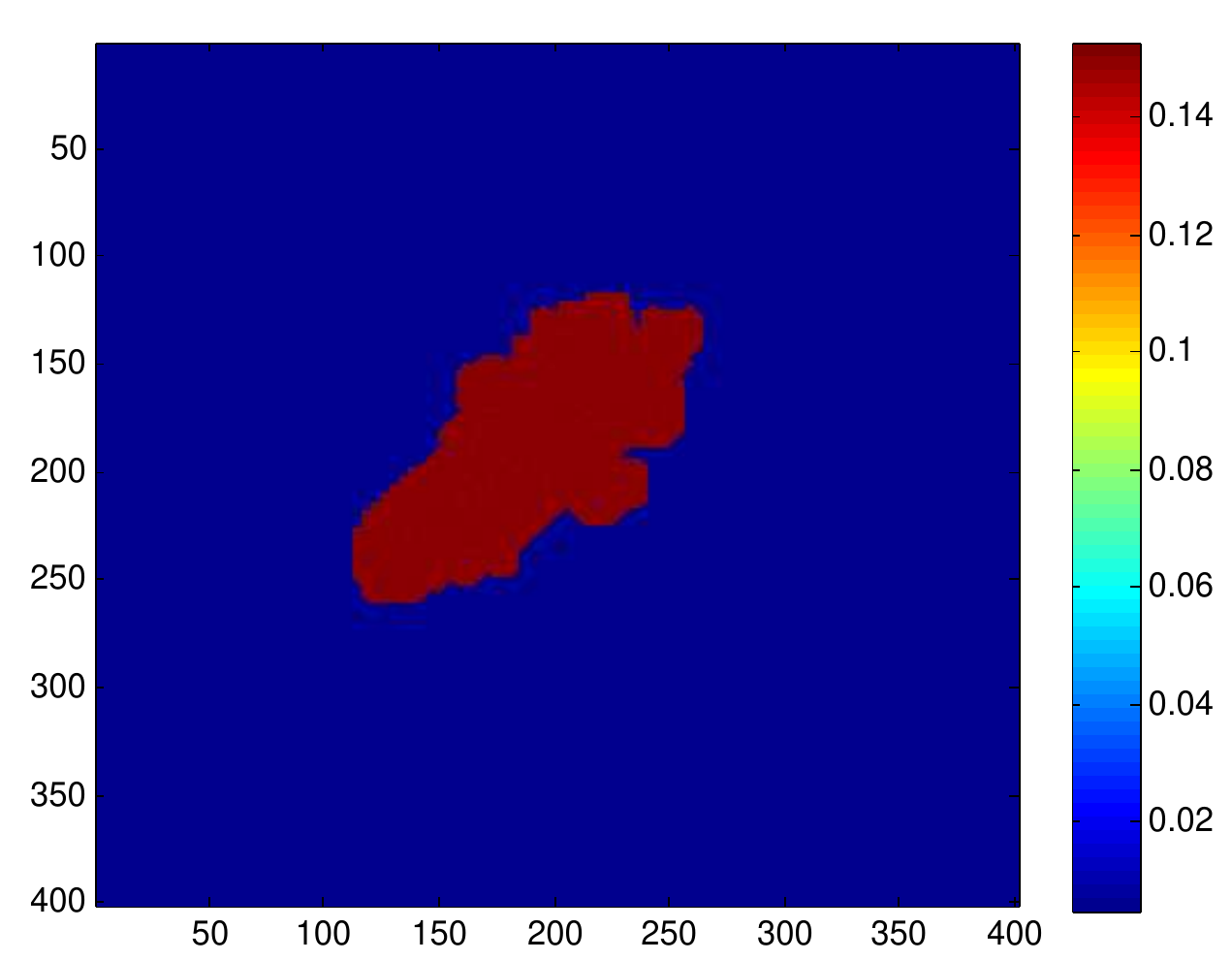}
\\
(a) Original shape.
\emp
\hspace{\stretch{1}}
\bmp{[t]}{1.5in}
\includegraphics[width=1.25in,height=1.25in]{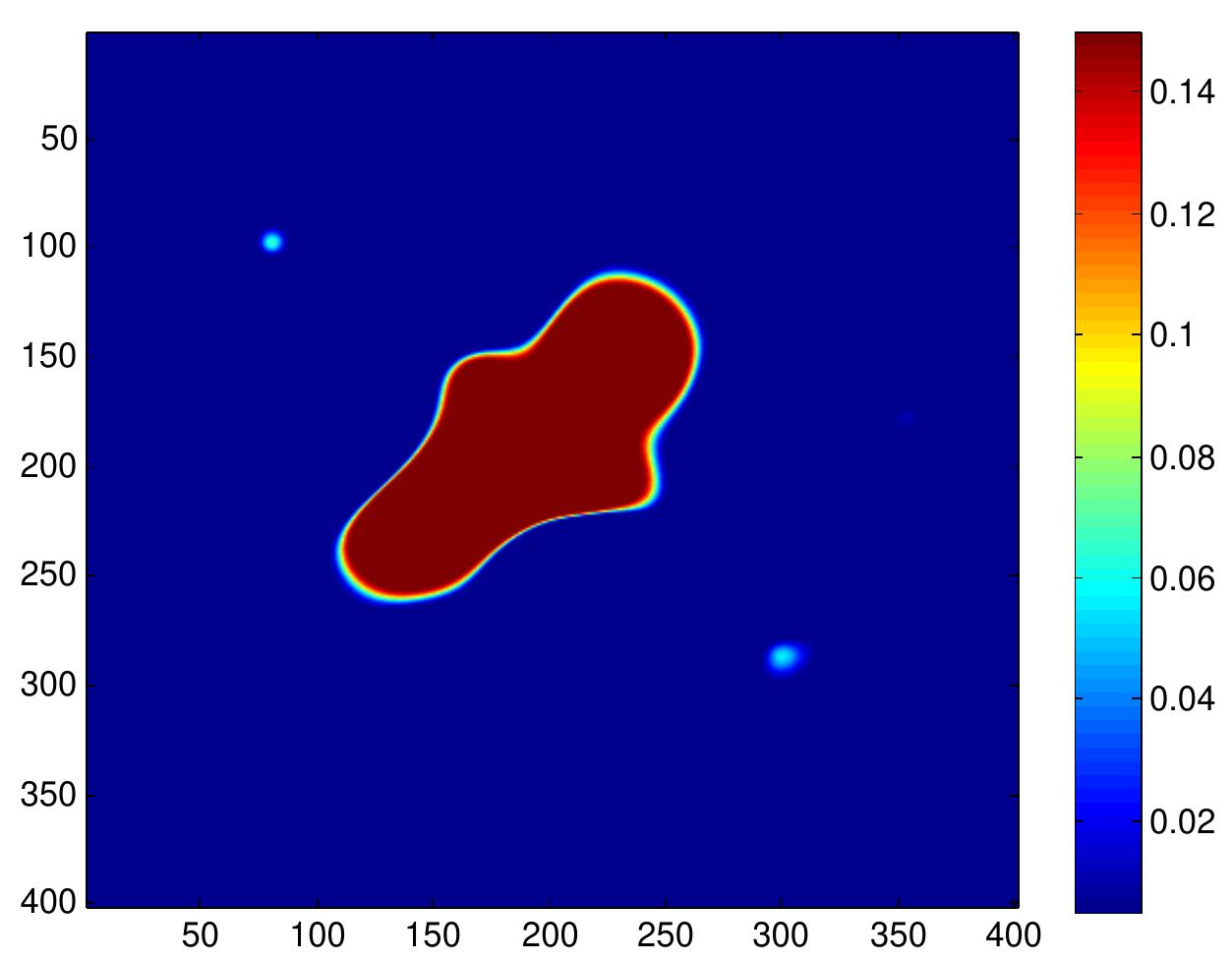}
\\
(b) Reconstruction using the exact objective function (full model).
\emp
\hspace{\stretch{1}}
\bmp{[t]}{1.5in}
\includegraphics[width=1.25in,height=1.25in]{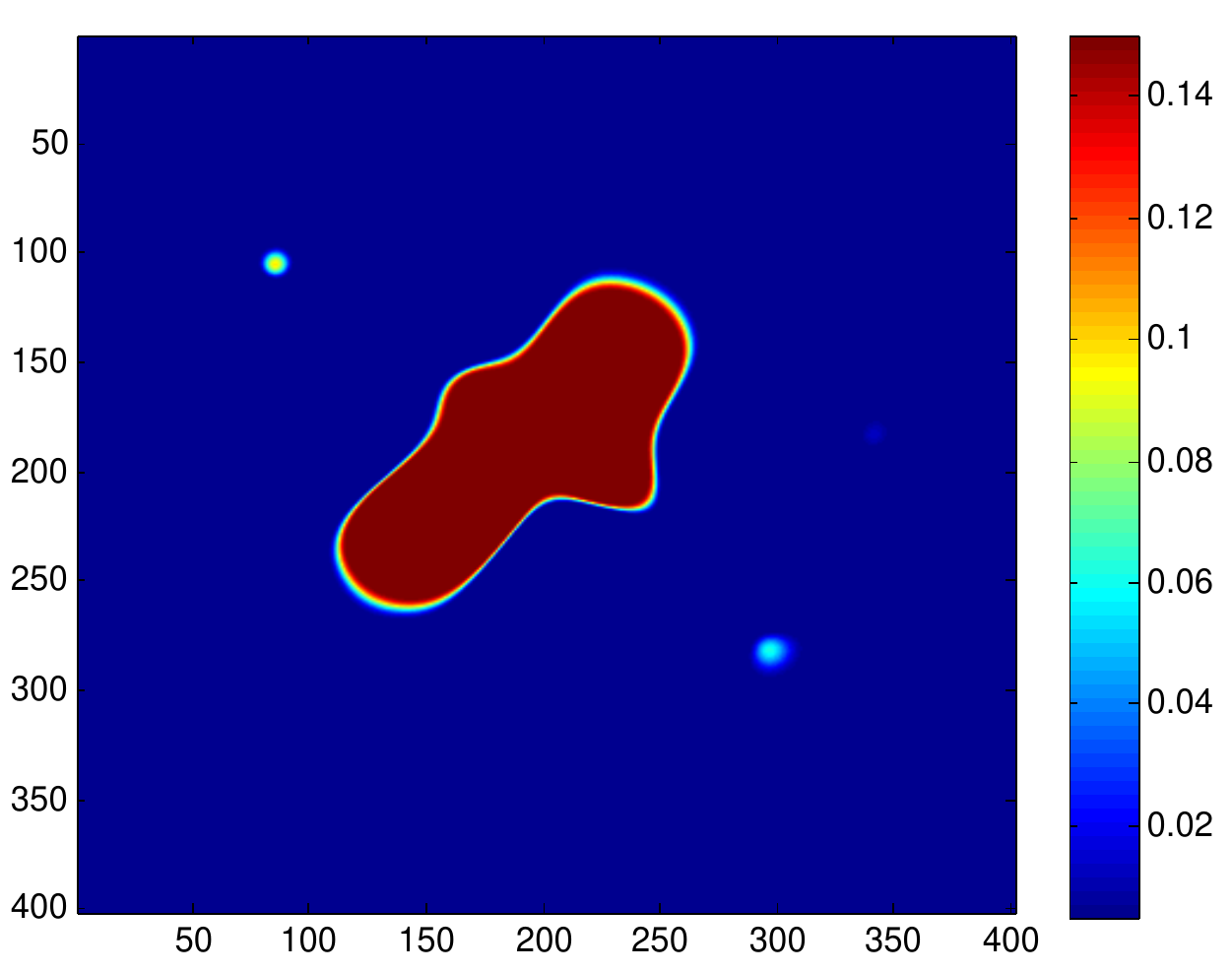}
\\
(c) Reconstruction using an approximate objective function (reduced order model).
\emp

\caption{Results for Problem~4. Reconstruction of a test anomaly on a $401 \times 401$ mesh, $160801$ degrees of freedom in the forward model, with $36$ sources and detectors, and $25$ basis functions. This reconstruction used the projection bases computed for the previous problem.
So, {\em only} the reduced model from the previous test case was used ($250$ degrees of freedom for the forward model).}
\label{fig:amoebe-recon}
\end{figure}

Next we consider two substantially larger problems. We now work on a $401 \times 401$ mesh resulting
in $160801$ degrees of freedom in the forward problem. The model in this case has $36$ sources and $36$ detectors.
We use $25$ CSRBFs for the parameterization
of the medium, which gives $100$ parameters for the nonlinear optimization problem.
The initial parameter vector is such that it produces the initial
absorption image given in Figure~\ref{fig:initial}. We use $5$ interpolation points in parameter space,
once again chosen from the first five iterations of the optimization using the full-order model.
Therefore, we solve $360$ large, sparse linear systems; one for each parameter point and each source as well as each detector. Combining both sets of solutions and using a rank revealing decomposition leads to a projection basis of $250$ vectors. So, the reduced model has order $250$,
and using the reduced model to approximate the objective function requires the solution of linear systems of size $250 \times 250$ rather than $160801 \times 160801$.
The original absorption image and the results of the reconstruction for Problem~3 are given in Figure~\ref{fig:cup-recon}. Solving the full-size forward problem, the exact objective function, for each step in the optimization requires the solution of $1120$ linear systems of dimension $160801$ versus only $360$ for the computing the reduced model. Then, using the reduced model for the forward problem requires only the solution of small linear systems of dimension $250$. The optimization using the reduced model takes $1216$ iterations. As the number of large linear systems is the main computational cost, we see that the computational cost solving a single image reconstruction problem is reduced by more than a factor three. Moreover, if multiple image reconstruction systems can be solved using the same reduced model, the computational savings are much higher.
The solution of Problem~4, which comes next, demonstrates that this is indeed the case.

For Problem~4, again we work with the same discretization of the forward problem, the same number of radial basis functions to model the absorption image, and the same number (and positions) of sources and detectors as for Problem~3. We also use the same initial guess for the parameter vector.
Figure~\ref{fig:amoebe-recon} gives (a) the original absorption image as well as (b) the reconstruction using the exact objective function (or full model) and (c) the reconstruction using the reduced model.
As before, we did not need to compute a new basis for the reduced model; we used the same basis as
computed for Problem~3.
The optimization using the full model required the solution of $896$ large linear systems of dimension $160801$, whereas the optimization using the reduced model required only the solution of $992$ small linear systems of dimension $250$. No large linear systems needed to be solved in the optimization using the reduced model.

Several types of scaling will determine the usefulness of our approach.
For this problem, we demonstrate one, the ratio between the 
number of degrees of freedom in the full model (exact objective function) 
and those in the reduced model (approximate objective function).
In table~\ref{tab:full-vs-reduced-dim}, we give the size of the reduced model required to get a similar error reduction in the image reconstruction. As the table shows, the dimension of the reduced model grows only modestly as a function of the dimension of full model. This comparison was carried out for a fixed number of sources and detectors.

\begin{table}
  \centering
  \begin{tabular}{||c|c||}
    \hline
    full model dimension (n) & reduced model dimension (r) \\
    \hline
    2500   & 80    \\
    10201  & 100   \\
    40401  & 150   \\
    160801 & 250   \\
    \hline
  \end{tabular}
  \caption{Table of reduced model dimensions versus full model dimensions for the forward model.}\label{tab:full-vs-reduced-dim}
\end{table}

\section{Conclusions and Future Work}\label{sec:Conc}
We have shown here how parametrized interpolatory  model reduction can significantly
reduce the cost of the inversion process in Diffuse Optical Tomography by reducing dramatically the
computational cost of forward problems.
The key observation is that function and Jacobian evaluations  arising in the inversion process
correspond to evaluations of a frequency response function and its gradient.
This motivated consideration of system-theoretic reduction methods and 
the use of parametrized interpolatory reduced models as surrogates for the full
forward model.
Four numerical examples
illustrate the efficiency of our approach.  These interpolatory reduced models were found to  reduce significantly the computational cost of the overall inversion process while producing negligible degradation in the quality of outcomes.  Section \ref{sec:anal_gb} offers some motivation and the outlines of an analysis that describes why these methods may be expected to be effective.

In this paper, we have used global reduction bases for constructing the reduced order models.  This approach is simple and effective. In our setting,
once the reduction bases $\bfV$ and $\bfW$ are constructed, the main online cost is the computation of the reduced matrix $\bfA_r(\bfsfp)$
using $\bfW^T \bfA(\bfsfp) \bfV$.  As we explained in \S \ref{ssec:GlobalBases}, due to the special structure of $\bfA(\bfsfp)$ in the DOT setting, this computation is still very cheap compared to solving large linear systems. 
To improve further the efficiency of our method (especially for larger 3D imaging problems),  we 
are  investigating incorporation of Discrete Empirical Interpolation Methods (DEIM) \cite{ChatSorDEIM_siam2010} and  Local-DEIM  \cite{Willcox2013}  
in the construction of affine approximations of $\bfA(\bfsfp)$. 
We are also investigating direct interpolation of local reduced-order quantities, as in
\cite{Panzer_etal2010,AmsallemFarhat2011,Degroote2010}. 

As an alternative to the full matrix interpolation of the frequency response function,
$\bfPsi(\omega;\bfsfp)$,  that is employed here, we are also investigating the use of tangential interpolation
\cite{BauB09,gallivan2005mrm} of 
$\bfPsi(\omega;\bfsfp)$ for cases with large numbers of detectors and sources.


\vspace{6pt}
\noindent

\bibliographystyle{abbrv}
\bibliography{npimr}

\begin{thebibliography}{10}

\bibitem{Aghasi_etal11}
A.~Aghasi, E.~Miller, and M.~E. Kilmer.
\newblock Parametric level set methods for inverse problems.
\newblock {\em SIAM Journal on Imaging Science}, 4:618--650, 2011.

\bibitem{Amsallem2008}
D.~Amsallem and C.~Farhat.
\newblock Interpolation method for the adaptation of reduced-order models to
  parameter changes and its application to aeroelasticity.
\newblock {\em AIAA Journal}, 46:1803--1813, July 2008.

\bibitem{AmsallemFarhat2011}
D.~Amsallem and C.~Farhat.
\newblock An online method for interpolating linear parametric reduced-order
  models.
\newblock {\em SIAM Journal on Scientific Computing}, 33(5):2169--2198, 2011.

\bibitem{Antil2011}
H.~Antil, M.~Heinkenschloss, and R.~H.~W. Hoppe.
\newblock Domain decomposition and balanced truncation model reduction for
  shape optimization of the {Stokes} system.
\newblock {\em Optimization Methods and Software}, 26(4--5):643--669, 2011.

\bibitem{Antil2012}
H.~Antil, M.~Heinkenschloss, R.~H.~W. Hoppe, C.~Linsenmann, and A.~Wixforth.
\newblock Reduced order modeling based shape optimization of surface acoustic
  wave driven microfluidic biochips.
\newblock {\em Mathematics and Computers in Simulation}, 82(10):1986--2003,
  2012.

\bibitem{antoulas2005approximation}
A.~Antoulas.
\newblock {\em {Approximation of Large-Scale Dynamical Systems (Advances in
  Design and Control)}}.
\newblock Society for Industrial and Applied Mathematics Philadelphia, PA, USA,
  2005.

\bibitem{Ant2010imr}
A.~Antoulas, C.~Beattie, and S.~Gugercin.
\newblock Interpolatory model reduction of large-scale dynamical systems.
\newblock In J.~Mohammadpour and K.~Grigoriadis, editors, {\em Efficient
  Modeling and Control of Large-Scale Systems}. Springer-Verlag, 2010.

\bibitem{Arian2002}
E.~Arian, M.~Fahl, and E.~Sachs.
\newblock Trust-region proper orthogonal decomposition models by optimization
  methods.
\newblock In {\em Proceedings of the 41st IEEE Conference on Decision and
  Control}, pages 3300--3305, Las Vegas, NV, 2002. IEEE.

\bibitem{arridge}
S.~R. Arridge.
\newblock Optical tomography in medical imaging.
\newblock {\em Inverse Problems}, Vol. 16:R41--R93, 1999.

\bibitem{BauB09}
U.~Baur and P.~Benner.
\newblock Model reduction for parametric systems using balanced truncation and
  interpolation.
\newblock {\em at--Automatisierungstechnik}, 57(8):411--420, 2009.

\bibitem{BauBBG09}
U.~Baur, P.~Benner, C.~Beattie, and S.~Gugercin.
\newblock Interpolatory projection methods for parameterized model reduction.
\newblock {\em SIAM Journal on Scientific Computing}, 33:2489--2518, 2011.

\bibitem{beattie2010isi}
C.~Beattie, S.~Gugercin, and S.~Wyatt.
\newblock Inexact solves in interpolatory model reduction.
\newblock {\em Linear Algebra and its Applications}, 2011.
\newblock Appeared on-line as doi:10.1016/j.laa.2011.07.015.

\bibitem{bond2005pmo}
B.~N. Bond and L.~Daniel.
\newblock {Parameterized model order reduction of nonlinear dynamical systems}.
\newblock In {\em IEEE/ACM Internat.Conf. on Computer-Aided Design, 2005.
  ICCAD-2005}, pages 487--494, 2005.

\bibitem{borcea2012model}
L.~Borcea, V.~Druskin, A.~V. Mamonov, and M.~Zaslavsky.
\newblock A model reduction approach to numerical inversion for a parabolic
  partial differential equation.
\newblock {\em arXiv preprint arXiv:1210.1257}, 2012.

\bibitem{BuiThanh2008}
T.~Bui-Thanh, K.~Willcox, and O.~Ghattas.
\newblock Model reduction for large-scale systems with high-dimensional
  parametric input space.
\newblock {\em SIAM Journal on Scientific Computing}, 30(6):3270--3288, 2008.

\bibitem{BuiThanh2008_AIAA}
T.~Bui-Thanh, K.~Willcox, and O.~Ghattas.
\newblock Parametric reduced-order models for probabilistic analysis of
  unsteady aerodynamic applications.
\newblock {\em AIAA Journal}, 46(10):2520--2529, 2008.

\bibitem{burger2005survey}
M.~Burger and S.~Osher.
\newblock {A survey on level set methods for inverse problems and optimal
  design}.
\newblock {\em European Journal of Applied Mathematics}, 16(02):263--301, 2005.

\bibitem{bushberg2003essential}
J.~Bushberg, J.~Seibert, E.~Leidholdt~Jr, J.~Boone, and E.~Goldschmidt~Jr.
\newblock {\em The essential physics of medical imaging}, volume~30.
\newblock Lippincott Williams {\&} Wilkens, 2003.

\bibitem{carrera2005inverse}
J.~Carrera, A.~Alcolea, A.~Medina, J.~Hidalgo, and L.~Slooten.
\newblock {Inverse problem in hydrogeology}.
\newblock {\em Hydrogeology Journal}, 13(1):206--222, 2005.

\bibitem{ChatSorDEIM_siam2010}
S.~Chaturantabut and D.~C. Sorensen.
\newblock Nonlinear model reduction via discrete empirical interpolation.
\newblock {\em SIAM Journal on Scientific Computing}, 32(5):2737--2764, 2010.

\bibitem{Daniel2004}
L.~Daniel, C.~Ong, S.~Low, K.~Lee, and J.~White.
\newblock A multiparameter moment matching model reduction approach for
  generating geometrically parameterized interconnect performance models.
\newblock {\em IEEE Transaction on Computer-Aided Design of Integrated Circuits
  and Systems}, 23(5):678--693, 2004.

\bibitem{StuKil11c}
E.~de~Sturler and M.~E. Kilmer.
\newblock A regularized {Gauss-Newton} trust region approach to imaging in
  diffuse optical tomography.
\newblock {\em SIAM Journal on Scientific Computing}, 33:3057 -- 3086, 2011.

\bibitem{Degroote2010}
J.~Degroote, J.~Vierendeels, and K.~Willcox.
\newblock Interpolation among reduced-order matrices to obtain parameterized
  models for design, optimization and probabilistic analysis.
\newblock {\em International Journal for Numerical Methods in Fluids},
  63(2):207--230, 2010.

\bibitem{dorn2006level}
O.~Dorn and D.~Lesselier.
\newblock {Level set methods for inverse scattering}.
\newblock {\em Inverse Problems}, 22:R67, 2006.

\bibitem{Druskin2011solution}
V.~Druskin, V.~Simoncini, and M.~Zaslavsky.
\newblock Solution of the time-domain inverse resistivity problem in the model
  reduction framework {Part I}. {O}ne-dimensional problem with {SISO} data.
\newblock {\em SIAM Journal on Scientific Computing}, 35(3):A1621--A1640, 2013.

\bibitem{FenB08}
L.~Feng and P.~Benner.
\newblock A robust algorithm for parametric model order reduction based on
  implicit moment matching.
\newblock {\em Proc. Appl. Math. Mech.}, 7:1021501--1021502, 2008.

\bibitem{gallivan2005mrm}
K.~Gallivan, A.~Vandendorpe, and P.~{Van~Dooren}.
\newblock {Model reduction of {MIMO} systems via tangential interpolation}.
\newblock {\em {SIAM} J. Matrix Anal. Appl.}, 26(2):328--349, 2005.

\bibitem{gugercin2008hmr}
S.~Gugercin, A.~Antoulas, and C.~Beattie.
\newblock $\mathcal{H}_2$ model reduction for large-scale linear dynamical
  systems.
\newblock {\em SIAM Journal on Matrix Analysis and Applications},
  30(2):609--638, 2008.

\bibitem{gunupudi2003ppt}
P.~Gunupudi, R.~Khazaka, M.~Nakhla, T.~Smy, and D.~Celo.
\newblock {Passive parameterized time-domain macromodels for high-speed
  transmission-line networks}.
\newblock {\em IEEE Trans. Microwave Theory and Techniques}, 51(12):2347--2354,
  2003.

\bibitem{haasdonk2011erm}
B.~Haasdonk and M.~Ohlberger.
\newblock Efficient reduced models and a posteriori error estimation for
  parametrized dynamical systems by offline/online decomposition.
\newblock {\em Mathematical and Computer Modelling of Dynamical Systems},
  17(2):145--161, 2011.

\bibitem{HabAsch00}
E.~Haber, U.~M. Ascher, and D.~Oldenburg.
\newblock On optimization techniques for solving nonlinear inverse problems.
\newblock {\em Inverse Problems}, 16:1263--1280, 2000.

\bibitem{Hay09}
A.~Har, J.~Borggaard, and D.~Pelletier.
\newblock Local improvements to reduced-order models using sensitivity analysis
  of the proper orthogonal decomposition.
\newblock {\em Journal of Fluid Mechanics}, 629:41--72, 2009.

\bibitem{hinze2005proper}
M.~Hinze and S.~Volkwein.
\newblock Proper orthogonal decomposition surrogate models for nonlinear
  dynamical systems: Error estimates and suboptimal control.
\newblock In {\em Dimension Reduction of Large-Scale Systems}, pages 261--306.
  Springer, 2005.

\bibitem{james33optimal}
A.~James, W.~Graham, K.~Hatfield, P.~Rao, and M.~Annable.
\newblock Optimal estimation of residual non--aqueous phase liquid saturations
  using partitioning tracer concentration data.
\newblock {\em Water Resources Research}, 33(12), 1997.

\bibitem{Kunisch2008}
K.~Kunisch and S.~Volkwein.
\newblock Proper orthogonal decomposition for optimality systems.
\newblock {\em ESAIM: Mathematical Modelling and Numerical Analysis},
  42(1):1--23, 1 2008.

\bibitem{liu2003computational}
G.~Liu and X.~Han.
\newblock {\em {Computational inverse techniques in nondestructive
  evaluation}}.
\newblock CRC, 2003.

\bibitem{louis1992medical}
A.~Louis.
\newblock {\em {Medical imaging: state of the art and future development}},
  volume~8.
\newblock Institute of Physics Publishing, 1992.

\bibitem{marklein2002linear}
R.~Marklein, K.~Mayer, R.~Hannemann, T.~Krylow, K.~Balasubramanian,
  K.~Langenberg, and V.~Schmitz.
\newblock {Linear and nonlinear inversion algorithms applied in nondestructive
  evaluation}.
\newblock {\em Inverse Problems}, 18:1733--1759, 2002.

\bibitem{nguyen2008best}
N.~Nguyen, A.~Patera, and J.~Peraire.
\newblock A Ôbest pointsÕ interpolation method for efficient approximation of
  parametrized functions.
\newblock {\em International Journal for Numerical Methods in Engineering},
  73(4):521--543, 2008.

\bibitem{osher1988fronts}
S.~Osher and J.~Sethian.
\newblock {Fronts propagating with curvature-dependent speed: algorithms based
  on Hamilton-Jacobi formulations}.
\newblock {\em Journal of computational physics}, 79(1):12--49, 1988.

\bibitem{Panzer_etal2010}
H.~Panzer, J.~Mohring, R.~Eid, and B.~Lohmann.
\newblock Parametric model order reduction by matrix interpolation.
\newblock {\em at--Automatisierungstechnik}, 58(8):475--484, 2010.

\bibitem{Willcox2013}
B.~Peherstorpher, D.~Butnaru, K.~Willcox, and H.-J. Bungartz.
\newblock Localized discrete emprical interpolation method.
\newblock Technical Report TR-13-1, MIT Aerospace Computational Design
  Laboratory, June 2013.

\bibitem{prud2002reliable}
C.~Prud'homme, D.~Rovas, K.~Veroy, L.~Machiels, Y.~Maday, A.~Patera, and
  G.~Turinici.
\newblock Reliable real-time solution of parametrized partial differential
  equations: Reduced-basis output bound methods.
\newblock {\em Journal of Fluids Engineering}, 124:70--80, 2002.

\bibitem{RozHP08}
G.~Rozza, D.~Huynh, and A.~Patera.
\newblock Reduced basis approximation and a posteriori error estimation for
  affinely parametrized elliptic coercive partial differential equations:
  application to transport and continuum mechanics.
\newblock {\em Arch. Comput. Methods Eng.}, 15(3):229--275, 2008.

\bibitem{santosa1996level}
F.~Santosa.
\newblock {A level-set approach for inverse problems involving obstacles}.
\newblock {\em ESAIM: Control, Optimisation and Calculus of Variations},
  1:17--33, 1996.

\bibitem{snieder1999inverse}
R.~Snieder and J.~Trampert.
\newblock {\em {Inverse problems in geophysics}}.
\newblock Springer, 1999.

\bibitem{stavroulakis2001inverse}
G.~Stavroulakis.
\newblock {\em {Inverse and crack identification problems in engineering
  mechanics}}.
\newblock Kluwer Academic Pub, 2001.

\bibitem{sun1994inverse}
N.~Z. Sun.
\newblock {\em Inverse problems in groundwater modeling}.
\newblock Kluwer Academic Publishers, Dordrecht, The Netherlands, 1994.

\bibitem{DoelAsch06}
K.~van~den Doel and U.~M. Ascher.
\newblock On level set regularization for highly ill-posed distributed
  parameter estimation problems.
\newblock {\em Journal of Computational Physics}, 216:707--723, 2006.

\bibitem{DoelAsch07}
K.~van~den Doel and U.~M. Ascher.
\newblock Dynamic level set regularization for large distributed parameter
  estimation problems.
\newblock {\em Inverse Problems}, 23:1271--1288, 2007.

\bibitem{DoelAsch10}
K.~van~den Doel, U.~M. Ascher, and A.~Leit{\~{a}}o.
\newblock Multiple level sets for piecewise constant surface reconstruction in
  highly ill-posed problems.
\newblock {\em Journal of Scientific Computing}, 43:44--66, 2010.

\bibitem{veroy2003posteriori}
K.~Veroy, C.~Prud'homme, D.~Rovas, and A.~Patera.
\newblock A posteriori error bounds for reduced-basis approximation of
  parametrized noncoercive and nonlinear elliptic partial differential
  equations.
\newblock In {\em Proceedings of the 16th AIAA Computational Fluid Dynamics
  Conference}, 2003.

\bibitem{Vogel}
C.~R. Vogel.
\newblock {\em Computational Methods for Inverse Problems}.
\newblock SIAM, Philadelphia, 2002.

\bibitem{webb2003introduction}
A.~Webb and G.~Kagadis.
\newblock {\em Introduction to biomedical imaging}, volume~30.
\newblock Wiley-IEEE Press, 2003.

\bibitem{wendland2005scattered}
H.~Wendland.
\newblock {\em {Scattered data approximation}}.
\newblock Cambridge University Press, 2005.

\bibitem{yeh22review}
W.~Yeh.
\newblock Review of parameter identification procedures in groundwater
  hydrology: The inverse problem.
\newblock {\em Water Resources Research}, 22(2), 1986.

\bibitem{yue2013}
Y.~Yue and K.~Meerbergen.
\newblock Accelerating optimization of parametric linear systems by model order
  reduction.
\newblock {\em SIAM Journal on Optimization}, 23(2):1344--1370, 2013.

\bibitem{zhdanov2002geophysical}
M.~Zhdanov.
\newblock {\em {Geophysical inverse theory and regularization problems}}.
\newblock Elsevier Science Ltd, 2002.

\end{thebibliography}

\end{document}